\documentclass[10pt]{article}

\usepackage{amsmath,amssymb}
\usepackage[a4paper]{geometry}

\newtheorem{Lemme}{Lemma}[section] 
\newtheorem{theorem}[Lemme]{Theorem}
\newtheorem{proposition}[Lemme]{Proposition}
\newtheorem{lemma}[Lemme]{Lemma}   

\newtheorem{remark}[Lemme]{Remark}	
\newtheorem{definition}[Lemme]{Definition}

\def\NN{{I~\hspace{-1.45ex}N} }

\def\Box{\leavevmode\vrule height 5pt width 4pt depth 0pt\relax}
\date{\today}

\begin{document}

\begin{center}
{\Large \textbf{Optimal control of steady second grade fluids with Dirichlet boundary conditions}}

\vspace{1 cm}
\textsc{Nadir Arada}{\footnote{Centro de Matem\'atica e Aplica\c{c}\~oes, Departamento de Matem\'atica, Faculdade de Ci\^encias e Tecnologia, Universidade Nova de Lisboa, Portugal. E-mail: naar{\char'100}fct.unl.pt}} 

%\vspace{1cm}

%\today

\end{center}

\vspace{ 1 cm}
%%%%%%%%%%%%%%%%%%%%%%%%%%%%%%%%%%%%%%%%%%%%%%%%%%%%%%%
\noindent \begin{abstract} \noindent  We consider optimal control  problems governed by systems describing the flow of an incompressible second grade fluid with Dirichlet boundary conditions. We prove the existence of an optimal solution, derive the corresponding  necessary optimality conditions and analyze its asymptotic behavior when the viscoelastic parameter tends to zero.   \vspace{2mm}\vspace{3mm}\vspace{3mm}
\end{abstract} 
%%%%%%%%%%%%%%%%%%%%%%%%%%%%%%%%%%%%%%%%%%%%%%%%%%%%%%%
\section{Introduction}\label{introduction}
\setcounter{equation}{0}
%%%%%%%%%%%%%%%%%%%%%%%%%%%%%%%%%%%%%%%%%%%%%%%%%%%%%%%
One of the important feature of complex non-Newtonian fluids is their ability to exhibit normal stress differences in simple shear flows, leading to characteristic phenomena like {\it rod-climbing} or {\it die-swell}. The second grade fluid model forms a subclass of differential type fluids of  complexity 2, and is one of the simplest constitutive models for flows of non-Newtonian fluids that can predict normal stress differences (cf. \cite{RE55} or \cite{NT65}). The corresponding stress is just a function of the pressure, the velocity gradient and some number of its higher material time derivatives (the Rivlin-Ericksen tensors). As a consequence, only an infinitesimal part of the history of the deformation gradient has an influence on the stress and, while they are good at predicting {\it creep}, these models cannot capture  {\it stress relaxation}. Nevertheless, due to their relative mathematical simplicity, there has been a great deal of interest on these models in recent years as they have been used successfully to predict slow steady motions of slurry flows, food rheology or flow of a water solution of polymers, where relaxation effects frequently seem to be rather insignificant. 
The corresponding equations of motion have the form
\begin{equation}\label{equation_sg}
	 \begin{array}{lll}\partial_t\left(\mathbf y-\alpha_1\Delta  \mathbf y\right)-\nu \Delta  \mathbf y\hspace{-3mm}&+
	\mathbf{curl}\left( \mathbf y-\left(2\alpha_1+\alpha_2\right)\Delta  \mathbf y
	\right)\times  \mathbf y\vspace{2mm}\\
	&+\left(\alpha_1+\alpha_2\right)\left(
	2\mathbf y\cdot \nabla \left(\Delta \mathbf y\right)-\Delta\left(\mathbf y\cdot \nabla \mathbf y\right)\right)
	+\nabla \pi= \mathbf u&  \mbox{in} \ \Omega,\end{array}
	\end{equation}
 where $\mathbf y$ is the velocity field, $\alpha_1$ and $\alpha_2$ are viscoelastic parameters (normal stress moduli), $\nu$ is the viscosity of the fluid,  $\pi$ is the hydrodynamic pressure, $\mathbf u$ is a given body force and $\Omega\subset \mathbb{R}^2$ is a bounded domain with boundary $\Gamma$. As this equation is set in dimension two, the vector $\mathbf y$ is written in the form $\mathbf y=(y\equiv (y_1,y_2),0)$ in order to define the vector product and the curl. Recall that in two dimensions, $ \mathrm{curl} \,  y=\tfrac{\partial y_2}{\partial x_1}-\tfrac{\partial y_1}{\partial x_2}$ and thus $\mathbf{curl} \, \mathbf y=(0,0, \mathrm{curl} \,  y)$. \vspace{1mm}\\
 According to \cite{DF74}, if the fluid modelled by equation (\ref{equation_sg}) is to be compatible with thermodynamics  in the sense that all motions of the fluid meet the Clausius-Duhem inequality and the assumption that the specific Helmholtz free energy of the fluid is a minimum in equilibrium, then 
	$$\nu\geq 0, \qquad
	\alpha_1\geq 0,\qquad \alpha_1+
	\alpha_2=0.$$
We refer to \cite{DR95} for a critical and extensive historical review of second-order fluid models and, in particular, for a discussion on the sign of the normal stress moduli. Here we will restraint to the simplified case $\alpha_1+\alpha_2=0$, with $\alpha_1\geq 0$ and $\nu>0$. Setting $\alpha_1=\alpha$, we can see that the problem of determining the velocity field $\mathbf{y}$ and the associated pressure $\pi$ satisfying the equations governing the flow of an incompressible second grade fluid reduces to 
$$
	\left\{ \begin{array}{ll}\partial_t\left(\mathbf y-\alpha\Delta  \mathbf y\right)-\nu \Delta  \mathbf y+
	\mathbf{curl}\left( \mathbf y-\alpha\Delta  \mathbf y
	\right)\times  \mathbf y+\nabla \pi= \mathbf u& \mbox{in} \ \Omega,\vspace{2mm} \\
             \mathrm{div} \,  \mathbf y=0& \mbox{in} \ \Omega.\end{array}\right.
	$$
In the inviscid case ($\nu=0$), the second-grade fluid equations are called $\alpha$-Euler equations. Initially proposed as a regularization of the incompressible Euler equations, they are geometrically significant and have been interpreted as a model of turbulence (cf. \cite{HMR981} and \cite{HMR98}). They also inspired another variant, called the $\alpha$-Navier-Stokes equations that turned out to be very relevant in turbulence modeling (cf. \cite{FHT1}, \cite{FHT2} and the references therein). These equations contain the regularizing term $-\nu\Delta\left(\mathbf y-\alpha \Delta\mathbf y\right)$ instead of $\nu\Delta\mathbf y$, making the dissipation stronger and the problem much easier to solve than in the case of second-grade fluids. When $\alpha=0$, the $\alpha$-Navier-Stokes and the second grade fluid equations are equivalent to the Navier-Stokes equation. \vspace{1mm}\\
Since the nonlinear term involves derivatives with higher order than the ones appearing in the viscous term, solving this problem is very challenging. The two dimensional case has been systematically studied for the first time in \cite{O81} and \cite{CO84} for both steady and unsteady cases with homogeneous Dirichlet boundary conditions. A Galerkin's method in the basis of the eigenfunctions of the operator $\mathbf{curl}(\mathbf{curl}(\mathbf y-\alpha\Delta \mathbf y))$ was especially designed to decompose the problem into a mixed elliptic-hyperbolic type, looking for the velocity $\mathbf y$ as a solution of a Stokes-like system coupled to a transport equation satisfied by $\mathbf{curl}\left( \mathbf y-\alpha\Delta  \mathbf y
	\right)$. Under minimal restrictions on the data, this approach allows the authors to establish the existence of solutions (and automatically recover  $H^3$ regularity) in the steady case, and to prove that the time-dependent version admits
 a unique global solution in the two dimensional case. This problem  received a lot of attention since these pioneering results and, without ambition for completeness, we refer to \cite{B99} where existence of a solution in the three dimensional steady case was established under a restriction on the size of the data. We also cite the extensions in \cite{GS94} and \cite{CG97}, where the three dimensional unsteady case was considered: global in time existence for small data was established, the former work using a Schauder fixed point argument while the latter considers the decomposition method on the system of Galerkin equations previously mentionned.\vspace{1mm}\\
This paper deals with the mathematical analysis of an optimal control problem associated with a steady viscous, incompressible second grade fluid. Control is effected through a distributed mechanical force and the objective is to match the velocity field to a given target field. More precisely, the controls and states are constrained to satisfy the following system of partial differential equations
	\begin{equation}\label{equation_etat}
	\left\{ \begin{array}{ll}-\nu \Delta  \mathbf y+
	\mathbf{curl}\left( \mathbf y-\alpha\Delta  \mathbf y
	\right)\times  \mathbf y+\nabla \pi= \mathbf u& \mbox{in} \ \Omega,\vspace{2mm} \\
             \mathrm{div} \,  \mathbf y=0& \mbox{in} \ \Omega,\vspace{2mm}\\
	\mathbf y=0&
	\mbox{on} \ \Gamma\end{array}\right.\end{equation}
and the optimal control problem reads as
	$$(P_\alpha) \
	 \left\{\begin{array}{ll}\mbox{minimize} & \displaystyle J(u,y)=
        \tfrac{1}{2}\int_\Omega\left|y-y_d\right|^2\,dx+\tfrac{\lambda}{2}
	\int_\Omega\left|u\right|^2\,dx\vspace{1mm}\\
	\mbox{subject to} & (u,y)\in U_{ad}\times H^3(\Omega) \ \mbox{satisfies} \
 (\ref{equation_etat}) \ \mbox{for some} \ \pi\in L^2(\Omega),\end{array}\right.$$
where  $\lambda\geq 0$, $y_d$ is some desired velocity field in $L^2(\Omega)$ and $U_{ad}$, the set of admissible controls, is a nonempty closed convex subset of 
$H(\mathrm{curl};\Omega)=\left\{v\in L^2(\Omega)\mid \mathrm{curl}\, u\in L^2(\Omega) \right\}$. \vspace{1mm}\\
Deriving the optimality conditions for problems governed by highly nonlinear equations is not an easy task  (cf. \cite{A12}, \cite{A13}, \cite{A14}, \cite{S05}, \cite{WR10}). The main difficulties are encountered when studying the solvability of the corresponding linearized and adjoint equations and are closely related with the regularity of the coefficients in the main part of the associated differential operators. As already mentioned, the choice of the special Galerkin basis used to study the state equation  is optimal in the sense that it allows us to prove the existence of regular solutions with minimal assumptions on the data. However,  the {\it direct} application of this approach to study the linearized and adjoint equations does not seem appropriate and does come at a cost. The main disadvantage is that it {\it automatically imposes} the derivation of a  $H^3$ estimate and this may be achieved only if high order derivatives of the state variable are well defined and if we impose additional restriction on their size. This in turn is only guaranteed if we consider regular, size constrained controls.  To overcome this difficulty, our idea is to consider an approximate optimal control problem governed by a state equation involving regularized controls. More precisely, if $(\bar u,\bar y)$ is a solution of $(P_\alpha)$ and $\varepsilon$ is a positive parameter, we consider the control problem
$$(P_\alpha^\varepsilon) \ \left\{\begin{array}{ll} \mbox{minimize} & I(u,y^\varepsilon)=\displaystyle J(u,y^\varepsilon)+\tfrac{1}{2}\int_\Omega
	\left|u-\bar u\right|^2\,dx+\tfrac{1}{2}\int_\Omega
	\left|{\rm curl}\left(u-\bar u\right)\right|^2\,dx
	\vspace{2mm}\\
	\mbox{subject to}  & (u,y^\varepsilon)\in U_{ad}\times 
	H^3(\Omega) \ \mbox{such that}
	\vspace{2mm}\\
&\left\{ \begin{array}{ll}-\nu \Delta  \mathbf y^\varepsilon+
	\mathbf{curl}\left(\mathbf y^\varepsilon-\alpha\, \Delta\mathbf y^\varepsilon\right)
		\times  \mathbf y^\varepsilon+\nabla \pi^\varepsilon= \varrho_\varepsilon(\mathbf u)& \mbox{in} \ \Omega,\vspace{2mm} \\
             \mathrm{div} \,  \mathbf y^\varepsilon=0& \mbox{in} \ \Omega,\vspace{2mm}\\
	\mathbf y^\varepsilon=0&
	\mbox{on} \ \Gamma,\end{array}\right.
\end{array}
\right.$$
where $\varrho_\varepsilon$ denotes a Friedrich mollifier. A careful analysis enables us to handle the issues mentioned above and to derive the corresponding approximate optimality conditions under natural restriction on the control variable. By passing to the limit in the regularization parameter $\varepsilon$, we recover the optimality conditions for $(P_\alpha)$.\vspace{1mm}\\
In this paper, we are also interested in the asymptotic behavior of the solutions of $(P_{\alpha})$, when the viscoelastic parameter $\alpha$ tends 
to zero. We will prove in particular that
	\begin{equation}\label{stability}\lim_{\alpha\rightarrow 0^+}\min(P_\alpha)=\min(P_0),\end{equation}
 where $(P_0)$ is the optimal control problem governed by the steady Navier-Stokes equations and defined by
	$$(P_0) \ \left\{\begin{array}{ll} \mbox{minimize} & J(u,y)
	\vspace{2mm}\\
	\mbox{subject to}  & (u,y)\in U_{ad}\times H^1(\Omega) \ \mbox{such that}
	\vspace{2mm}\\
&\left\{\begin{array}{ll}-\nu \Delta y+
	y\cdot \nabla y+\nabla \pi=u& \quad 
	\mbox{in} \ \Omega,\vspace{2mm}\\
	\mathrm{div}\, y=0& \quad 
	\mbox{in} \ \Omega,\vspace{2mm}\\
	y=0& \quad 
	\mbox{on} \ \Gamma.\end{array}\right.\end{array}
\right.$$
To obtain such a result, we first establish that the sequence of solutions $(y_\alpha)_\alpha$  of (\ref{equation_etat}) converges to $y$, a solution of the Navier-Stokes equation, when $\alpha$ tends to zero. Next we prove that if $(\bar u_\alpha,\bar y_\alpha)$ is a solution to the problem $(P_\alpha)$ then the sequence $(\bar u_\alpha,\bar y_\alpha)_\alpha$ converges  to a solution $(\bar u_0,\bar y_0)$ of $(P_0)$. Another aspect concerns the necessary optimality conditions. To study the asymptotic behavior of these conditions, we analyze the adjoint equations for $(P_\alpha)$ and prove that the sequence of adjoint solutions converges to the solution of the adjoint equation for $(P_0)$. The optimality conditions for $(P_0)$ are then obtained by passing to the limit in the optimality conditions for $(P_\alpha)$\vspace{1mm}\\
The plan of the present paper is as follows. The main results are stated in Section 2. Notation and preliminary results related with the nonlinear terms are given in Section 3. Section 4 is devoted to the existence and uniqueness results for the state and the linearized state equation and to the derivation of the corresponding estimates. In Section 5, we analyze the Lipschitz continuity and the G\^ateaux differentiability of  the control-to-state mapping and we consider the solvability of the adjoint equation in Section 6. Finally, the proof of the main results are given in Section 7. 

%%%%%%%%%%%%%%%%%%%%%%%%%%%%%%%%%%%%%%%%%
\section{Statement of the main results}
\setcounter{equation}{0}
%%%%%%%%%%%%%%%%%%%%%%%%%%%%%%%%%%%%%%%%%
We first establish the existence of optimal solutions for problem $(P_\alpha)$. 
%%%%%%%%%%%%%%%%%%%%%%%%%%%%%%%%%%%%%%%%%
\begin{theorem} \label{main_existence}Assume that $U_{ad}$ is bounded in $H(\mathrm{curl};\Omega)$.
Then problem $(P_\alpha)$ admits at least a solution.
\end{theorem}
%%%%%%%%%%%%%%%%%%%%%%%%%%%%%%%%%%%%%%%%%
To derive the corresponding necessary optimality conditions (stated in the next result), we need to restrain the optimal control size. Such a restriction, well known and widely used when dealing with optimal control problems governed by the steady Navier-Stokes equations, should be set within the natural functional framework of $H({\rm curl};\Omega)$, without requiring additional regularity on the control. Besides the difficulties inherent to the highly nonlinear nature of the state equation, and its implications on the linearized and adjoint equations, this is one of the main issues we must overcome.
%%%%%%%%%%%%%%%%%%%%%%%%%%%%%%%%%%%%%%%%%
\begin{theorem} \label{main_1} Let $(\bar{ u}_\alpha,\bar{y}_\alpha)$ be a solution of $(P_\alpha)$. There exists a positive constant 
$\bar\kappa$, depending only on $\Omega$, such that if the following condition holds
	\begin{equation}\label{control_constraint}
	\bar\kappa \left(\left\|\bar u_\alpha\right\|_2+\alpha
	\left\|\mathrm{curl}\, \bar u_\alpha\right\|_2\right)<\nu^2\end{equation}
then there exists $\bar{p}_\alpha\in H^1(\Omega)$ weak solution of 
	\begin{equation}\label{adj_opt_eq_alpha}
	\left\{ \begin{array}{ll}-\nu \Delta  \bar{\mathbf p}_\alpha-
	\mathbf{curl}\,\sigma\left(\bar{\mathbf y}_\alpha\right)\times
	\bar{\mathbf p}_\alpha+\mathbf{curl}
	\left(\sigma\left(\bar{\mathbf y}_\alpha\times \bar{\mathbf p}_\alpha\right)\right)+\nabla \pi= \bar{\mathbf y}_\alpha-\mathbf y_d& \mbox{in} \ \Omega,\vspace{2mm} \\
             \mathrm{div} \,  \bar{\mathbf p}_\alpha=0& \mbox{in} \ \Omega,\vspace{2mm}\\
                \bar{\mathbf p}_\alpha=0& \mbox{on} \ \Gamma,\end{array}\right.\end{equation}
and satisfying
	\begin{equation}\label{opt_control_alpha}\left(\bar{ p}_\alpha+\lambda\bar{ u}_\alpha,v-\bar{ u}_\alpha\right)\geq 0 \qquad 
	\mbox{for all} \ v\in U_{ad}.\end{equation}
\end{theorem}
%%%%%%%%%%%%%%%%%%%%%%%%%%%%%%%%%%%%%%%%%%%%%%
Finally, we consider the asymptotic analysis of the optimal control $(P_\alpha)$. We first prove that if $(\bar u_\alpha,\bar y_\alpha)$ is a solution of $(P_\alpha)$, then a cluster point (for an appropriate topology) is a solution of problem $(P_0)$ and the stability property (\ref{stability}) holds. Moreorer, if $(\bar u_\alpha,\bar y_\alpha,\bar p_\alpha)$ is defined as in Theorem $\ref{main_1}$, then $\bar p_\alpha$ converges to some $\bar p_0$ satisfying the optimality conditions of problem $(P_0)$. More precisely, we have the following result.
%%%%%%%%%%%%%%%%%%%%%%%%%%%%%%%%%%%%%%%%%%%%%%
\begin{theorem} \label{assympt_1} Let $(\bar u_\alpha,\bar y_\alpha,\bar p_\alpha)$ defined as in Theorem $\ref{main_1}$. Then\vspace{2mm}\\
$i)$ $(\bar u_\alpha,\bar y_\alpha)$ strongly converges in $L^2(\Omega)\times V$ $($up to a subsequence when $\alpha$ tends to zero$)$ to a limit point $(\bar u_0,\bar y_0)$ solution o $(P_0)$.\vspace{2mm}\\
$ii)$ $\bar p_\alpha$ converges $($up to a subsequence when $\alpha$ tends to zero$)$ for the weak topology of $V$ to $\bar p_0$, weak solution of the adjoint equation
	\begin{equation}\label{adjoint_limit}
	\left\{\begin{array}{ll}-\nu \Delta \bar p_0-\bar y_0\cdot 
	\nabla \bar p_0+\left(\nabla \bar y_0\right)^\top \bar p_0+\nabla \pi=\bar y_0-y_d& \quad 
	\mbox{in} \ \Omega,\vspace{2mm}\\
	\mathrm{div}\, \bar p_0=0& \quad 
	\mbox{in} \ \Omega,\vspace{2mm}\\
	\bar p_0=0& \quad 
	\mbox{on} \ \Gamma\end{array}\right.
	\end{equation}
and satisfying the optimality condition
	$$\left(\bar p_0+\lambda\bar u_0,v-\bar u_0\right)\geq 0 \qquad 
	\mbox{for all} \ v\in U_{ad}.$$
\end{theorem}
%%%%%%%%%%%%%%%%%%%%%%%%%%%%%%%%%%%%%%%%%%%%%%
\begin{remark} Taking into account Remark $\ref{remark1}$ below, it follows that if $\bar u_\alpha$ satisfies $(\ref{control_constraint})$, then the limit $\bar u_0$ satisfies 
	$S_4^2 S_2\left\|\bar u_0\right\|_2<\nu^2$,
which implies the uniqueness of $\bar y_0$ and $\bar p_0$.\end{remark}
%%%%%%%%%%%%%%%%%%%%%%%%%%%%%%%%%%%%%%%%%%%%%%
\begin{remark} Unlike the two dimensional case where existence of at least a weak solution for the state equation can be established without restriction on the size of the data, the existence of such a solution is only guaranteed for small data in the three dimensional case $($see e.g. $\cite{B99}$$)$. As a consequence, the results stated in Theorem $\ref{main_existence}$, Theorem $\ref{main_1}$ and Theorem $\ref{assympt_1}$ may be extended to the three dimensional case, under additional restrictions on the whole set of admissible controls. \end{remark}
%%%%%%%%%%%%%%%%%%%%%%%%%%%%%%%%%%%%%%%%%%%%%%%%%%%%%%%
\section{Notation, assumptions and preliminary results}
\setcounter{equation}{0}
%%%%%%%%%%%%%%%%%%%%%%%%%%%%%%%%%%%%%%%%%%%%%%%%%%%%%%%
\subsection{Functional setting}
%%%%%%%%%%%%%%%%%%%%%%%%%%%%%%%%%%%%%%%%%%%%%%%%%%%%%%%%%%%%%
Throughout the paper $\Omega$ is a bounded, simply connected domain in $\mathbb{R}^2$. The boundary of $\Omega$ is denoted by $\Gamma$ and is of class $C^{2,1}$. The standard Sobolev spaces are denoted by $W^{k,p}(\Omega)$ ($k\in \mathbb{N}$ and $1<p<\infty$), and their norms by $\|\cdot\|_{k,p}$. We set $W^{k,2}(\Omega)\equiv H^k(\Omega)$ and $\|\cdot\|_{k,2}\equiv \|\cdot\|_{H^k}$. In order to simplify the presentation, we will use the notation 
	$$\sigma(v)=v-\alpha \Delta v, \qquad v\in H^2(\Omega)$$
in all the sequel. We will also frequently use the scalar product in $L^2(\Omega)$
	$$\left(u,v\right)=\displaystyle\int_\Omega u(x)\cdot v(x)\,dx,$$
the semi-norm of $H^1(\Omega)$
	$$\left|v\right|_{H^1}=\left\|\nabla v\right\|_2$$
and in order to eliminate the pressure in the different variational formulations, we will work in divergence-free spaces and consider the following Hilbert space
$$V=\left\{v\in H^1_0(\Omega)\mid \mathrm{div} \, v=0 \
 \mbox{in} \, \Omega\right\}.$$
Recall also the Poincar\'e and Sobolev inequalities, respectively given by
	$$\left\|v\right\|_2\leq S_{2} \left|v\right|_{H^1}
	\qquad \mbox{for all} \ v\in V,$$
	$$\left\|v\right\|_4\leq S_{4} \left|v\right|_{H^1}
	\qquad \mbox{for all} \ y\in V.$$
We introduce the space
	$$V_2=\left\{v\in V\mid {\rm curl}\, \sigma(v)\in L^2(\Omega)\right\}$$
equipped with the scalar product 
	$$\left(u,v\right)_{V_2}=\left(u,v\right)+\alpha \left(\nabla u,\nabla v\right)+\left({\rm curl}\, \sigma(u),{\rm curl}\, \sigma(v)\right)$$
and associated semi-norm
	$$\left|v\right|_{V_ 2}=\left\|{\rm curl}\, \sigma(v)\right\|_2.$$
We finally introduce the space (of controls) 
$$H({\rm curl};\Omega)=\left\{v\in L^2(\Omega)\mid {\rm curl} \, v\in L^2(\Omega)\right\}$$
equipped with the scalar product
	$$\left(u,v\right)_{H({\rm curl};\Omega)}=
	\left(u,v\right)+\left({\rm curl} \, u,{\rm curl} \, v\right)$$
and which is a Hilbert space for the associated norm
	$$\left\|v\right\|_{H({\rm curl};\Omega)}=\left(v,v\right)_{H({\rm curl};\Omega)}^{\frac{1}{2}}.$$
%%%%%%%%%%%%%%%%%%%%%%%%%%%%%%%%%%%%%%%%%%%%%%%%%%%%%%%
\subsection{Auxiliary results}
%%%%%%%%%%%%%%%%%%%%%%%%%%%%%%%%%%%%%%%%%%%%%%%%%%%%%%%
The aim of this section is to present some results that will be used throughout the paper. We first recall that the space $V_ 2$, particularly well adapted to handle the partial differential equations we are considering, is continuously embedded in $H^3(\Omega)$ (see e.g. \cite{CG97}).
\begin{lemma} Any $y\in V_2$ belongs to $H^3(\Omega)$ and there exists a constant $c(\alpha)$ such that
	$$\left\|y\right\|_{H^3}\leq c(\alpha) \left|y\right|_{V_2}.$$ 
\end{lemma}
The second lemma will be useful when dealing with a priori estimates for the linearized state and adjoint state equations.
%%%%%%%%%%%%%%%%%%%%%%%%%%%%%%%%%%%%%%%%%%%%%%%%%%%%%%%
\begin{lemma} \label{infty_V2} Let $y\in V_2$. Then, the following estimate holds
	$$\left\|y\right\|_{\infty}\leq \tfrac{c}{\alpha^\frac{1}{3}}
	\left|y\right|_{H^1}^{\frac{2}{3}}
	\left|y\right|_{V_2}^{\frac{1}{3}},$$
where $c$ is a positive constant only depending on $\Omega$.\end{lemma}
%%%%%%%%%%%%%%%%%%%%%%%%%%%%%%%%%%%%%%%%%%%%%%%%%%%%%%%
{\bf Proof.} Since $\mathrm{curl}\,\sigma(y)\in L^2(\Omega)$
and $\nabla\cdot \left(\mathrm{curl}\,\sigma(y)\right)=0$, there exists a unique vector-potential $\psi\in H^1(\Omega)$ such that
	$$\left\{\begin{array}{ll}\rm{curl}\,
	 \psi=\mathrm{curl}\,\sigma(y) &\quad \mbox{in} \ \Omega, \vspace{2mm}\\
	\nabla \cdot \psi=0 & \quad\mbox{in} \ \Omega,\vspace{2mm}\\
	\psi\cdot n=0 & \quad\mbox{on} \ \Gamma\end{array}\right.$$
and 
	\begin{equation}\label{sigma_phi}
	\left\|\psi\right\|_{H^1}\leq 
	c\left|y\right|_{V_2}.\end{equation}
It follows that
	$${\rm curl} \left(y-\alpha 
	\Delta y-\psi\right)=0$$
and the fact that $\Omega$ is simply connected  implies that
there exists $\pi\in L^2(\Omega)$ such that
	$$y-\alpha 
	\Delta y-\psi+\nabla \pi=0.$$
(For the proof of such a result, see Theorem 2.9, Chapter 1 in \cite{GR86}.)
Hence $y$ is the solution of the Stokes system
	$$-\Delta y+\nabla\left(\tfrac{\pi}{\alpha}\right)=\tfrac{1}{\alpha}
	\left(\psi-y\right)$$
and satisfies
	\begin{equation}\label{y_phi}\left\|y\right\|_{H^2}
	\leq \tfrac{c}{\alpha}
	\left\|\psi-y\right\|_2.\end{equation}
Observing that
	$$\left(\psi,y\right)=
	\left(y-\alpha \Delta y+\nabla \pi,y\right)=
	\left\|y\right\|_2^2+
	\alpha \left|y\right|_{H^1}^2,$$
we obtain
	$$\left\|\psi-y\right\|_2^2=
	\left\|\psi\right\|_2^2-
	\left\|y\right\|_2^2-2\alpha
	 \left|y\right|_{H^1}^2\leq
	 \left\|\psi\right\|_2^2.$$
Combining (\ref{sigma_phi}) and (\ref{y_phi}), we deduce that
\begin{equation}\label{yh_curl_sigma} \left\|y\right\|_{H^2}\leq 
	\tfrac{c}{\alpha}\left|y\right|_{V_2}.
	\end{equation}
Finally, the interpolation inequalities
	$$\left\|y\right\|_\infty\leq c 
	\left\|y\right\|_2^{\frac{1}{3}} 
	\left\|y\right\|_{1,4}^{\frac{2}{3}}\leq c 
	\left\|y\right\|_2^{\frac{1}{3}} 
	\left(\left\|y\right\|_{H^1}^\frac{1}{2}\left\|y\right\|_{H^2}^\frac{1}{2}\right)^{\frac{2}{3}},$$
together with the Poincar\'e inequality and (\ref{yh_curl_sigma}) yield
	$$\left\|y\right\|_\infty\leq 
	c\left|y\right|_{H^1}^{\frac{2}{3}}\, 
	\left\|y\right\|_{H^2}^{\frac{1}{3}}
	\leq \tfrac{c}{\alpha^\frac{1}{3}} 
	\left|y\right|_{H^1}^{\frac{2}{3}}\, 
	\left|y\right|_{V_2}^{\frac{1}{3}}$$
and the claimed result is proved. $\hfill \Box$\vspace{2mm}\\
%%%%%%%%%%%%%%%%%%%%%%%%%%%%%%%%%%%%%%%%%%%%%%%%%%%%%%%
%%%%%%%%%%%%%%%%%%%%%%%%%%%%%%%%%%%%%%%%%%%%%%%%%%%%%%%
The first identity in the next result is standard and relates the nonlinear term in (\ref{equation_etat}), and similar terms appearing in the  linearized and adjoint state equations, to the classical trilinear form used in the Euler and Navier-Stokes equations and defined by
	$$b(\phi,z,y)=\left(\phi\cdot \nabla z,y\right).$$  
The second identity deals with another term only appearing in the adjoint state equation. 
%%%%%%%%%%%%%%%%%%%%%%%%%%%%%%%%%%%%%%%%%%%%%%%%%%%%%%%
\begin{lemma} \label{non_lin_curl}
Let $y,z \in  V_2$ and $\phi\in  V$. Then
	$$\left( \mathbf{curl}\, \sigma(\mathbf y)
	\times  \mathbf z, \boldsymbol\phi\right)=
	b\left(\phi,z, \sigma(y)\right) 
	-b\left(z,\phi,\sigma(y)\right).$$
Let $y, z$ and $\phi$ be in $V_2$. Then 
	$$\left( \mathbf{curl}\, \sigma\left(\mathbf y\times \mathbf z\right), \boldsymbol\phi\right)
	=b\left(z,y, \sigma(\phi)\right)-b\left(y,z,\sigma(\phi)\right).$$
\end{lemma}
%%%%%%%%%%%%%%%%%%%%%%%%%%%%%%%%%%%%%%%%%%%%%%%%%%%%%%%
{\bf Proof.} By using a standard integration by parts, we can easily prove that for every $ y, z\in  V_2$
and every $\phi\in V$, we have
	$$\begin{array}{ll}\left( \mathbf{curl}\, \sigma(\mathbf y)
	\times  \mathbf z, \boldsymbol\phi\right)&=\left( \mathbf{curl}\, \sigma(\mathbf y)
	, \mathbf z\times \boldsymbol\phi\right)=\left( \sigma(\mathbf y)
	, \mathbf{curl}\left(\mathbf z\times \boldsymbol\phi\right)\right)\vspace{2mm}\\
	&=\left( \sigma(\mathbf y)
	,\boldsymbol\phi\cdot \nabla \mathbf z-\mathbf z\cdot \nabla \boldsymbol\phi\right)=b\left(\phi,z, \sigma(y)\right) 
	-b\left(z,\phi,\sigma(y)\right).\end{array}$$
and the first identity is proved. Similarly, for $y, z$ and $\phi$ be in $V_2$ we have
$$\begin{array}{ll}\left( \mathbf{curl}\, \sigma\left(\mathbf y\times \mathbf z\right), \boldsymbol\phi\right)&=
	\left( \mathbf{curl}\, \left(\mathbf y\times \mathbf z\right), \boldsymbol\phi\right)-\alpha
	\left( \mathbf{curl}\, \Delta\left(\mathbf y\times \mathbf z\right), \boldsymbol\phi\right)\vspace{2mm}\\
&=\displaystyle \left(\mathbf z\cdot \nabla \mathbf y-\mathbf y\cdot \nabla \mathbf z,\boldsymbol\phi\right)-\alpha\left(\Delta\left(\mathbf y\times \mathbf z\right),\mathbf{curl}\, \boldsymbol\phi\right)\vspace{2mm}\\
	&=b\left(z,y,\phi\right)-b\left(y,z,\phi\right)+
	\alpha\left(\mathbf{curl}\left(\mathbf{curl}\left(\mathbf y\times \mathbf z\right)\right)
	-\nabla\left( \mathrm{div}\left(\mathbf y\times \mathbf z\right)\right),\mathbf{curl}\, \boldsymbol\phi\right)\vspace{2mm}\\
	&=b\left(z,y,\phi\right)-b\left(y,z,\phi\right)+\alpha\left(\mathbf{curl}\left(\mathbf{curl}\left(\mathbf y\times \mathbf z\right)\right),\mathbf{curl}\,
	 \boldsymbol\phi\right)\vspace{2mm}\\
	&=b\left(z,y,\phi\right)-b\left(y,z,\phi\right)+\alpha \left(\mathbf{curl}\left(\mathbf y\times \mathbf z\right),\mathbf{curl}\left(\mathbf{curl}\, \boldsymbol\phi\right)\right)\vspace{2mm}\\
	&=b\left(z,y,\phi\right)-b\left(y,z,\phi\right)-\alpha\left(\mathbf{curl}\left(\mathbf y\times \mathbf z\right),\Delta \boldsymbol\phi-\nabla\left( \mathrm{div}\, \boldsymbol\phi\right)\right)\vspace{2mm}\\
	&=b\left(z,y,\phi\right)-b\left(y,z,\phi\right)-\alpha\left(b(\mathbf z,\mathbf y,\Delta\boldsymbol\phi)-b(\mathbf y,\mathbf z,\Delta \boldsymbol\phi)\right)\vspace{2mm}\\
	&=b\left(z,y, \sigma(\phi)\right)-b\left(y,z,\sigma(\phi)\right)
	\end{array}$$
and the second identity is proved.$\hfill \Box$\vspace{2mm}\\
%%%%%%%%%%%%%%%%%%%%%%%%%%%%%%%%%%%%%%%%%%%%%%%%%%%%%%%%%%%%
As will be seen in the sequel, the first identity in Lemma \ref{non_lin_curl} enables us to give an adequate variational setting for the state and linearized state equations. Based on the corresponding definitions, we can derive $H^1$ and $H^3$ a priori estimates and establish existence results. 
Similarly, combining the two identities in Lemma \ref{non_lin_curl}, we can propose a variational formulation for the adjoint equation and establish a $H^1$ estimate of the corresponding solution.
%%%%%%%%%%%%%%%%%%%%%%%%%%%%%%%%%%%%%%%%%%%%%%%%%%%%%%%
This section concludes with a result that will be used to establish a uniqueness
result for the state equation and to derive $H^1$ a priori estimates for the linearized state equation and the adjoint equation. 
%%%%%%%%%%%%%%%%%%%%%%%%%%%%%%%%%%%%%%%%%%%%%%%%%%%%%%%
%%%%%%%%%%%%%%%%%%%%%%%%%%%%%%%%%%%%%%%%%%%%%%%%%%%%%%%
\begin{lemma}\label{rm2}
Let $ y,z \in  V_2$. Then
	$$
	\left|\left( \mathbf{curl}\,  \sigma(\mathbf z)	\times  \mathbf y, \mathbf z\right)\right|
	\leq \left( S_4^2\left|y\right|_{H^1}+ \kappa  
	\alpha\left\|y\right\|_{H^3}\right)\left|z\right|_{H^1}^2,$$
where  $\kappa$ is a positive constant only depending on $\Omega$. 
\end{lemma}
%%%%%%%%%%%%%%%%%%%%%%%%%%%%%%%%%%%%%%%%%%%%%%%%%%%%%%%
{\bf Proof.} Lemma \ref{non_lin_curl} together with classical arguments show that
	$$\begin{array}{ll}
	\left( \mathbf{curl}\,\sigma(\mathbf z)\times  
	\mathbf y,\mathbf z\right)&=b\left(z, y, \sigma(z)\right) 
	-b\left(y,z, \sigma(z)\right) \vspace{2mm}\\
	&=b(z,y,z)-b(y,z,z)-\alpha 
	\left(z\cdot \nabla y-y\cdot \nabla z,\Delta z\right)\vspace{2mm}\\
	&=b(z,y,z)
	+\alpha\left(b\left(\mathbf z, \mathbf{curl}\, \mathbf y,\mathbf{curl}\, \mathbf z\right)
	 -b\left(\mathbf y, \mathbf{curl}\, \mathbf z,\mathbf{curl}\, \mathbf z\right)\right)\vspace{2mm}\\
	&\displaystyle +2\alpha\sum_{k=1}^3\left(\nabla \mathbf z_k\times \nabla \mathbf y_k,\mathbf{curl}\,\mathbf z\right)\vspace{-1mm}\\
	&\displaystyle =b(z,y,z)
	+\alpha b\left(\mathbf z, \mathbf{curl}\, \mathbf y,\mathbf{curl}\, \mathbf z\right)
	+2\alpha\sum_{k=1}^2\left(\nabla \mathbf z_k\times \nabla \mathbf y_k,\mathbf{curl}\,\mathbf z\right).
	\end{array}$$
Therefore,
	$$\begin{array}{ll}&\left|\left( \mathbf{curl}\,\sigma(\mathbf z)\times  \mathbf y,\mathbf z\right)\right|
	\vspace{2mm}\\
	&\displaystyle\leq \|z\|_4^2 \left|y\right|_{H^1}+
	\alpha\left(\|z\|_4 \left\|\nabla \mathbf{curl}\,\mathbf y\right\|_{4} \left\|\mathbf{curl}\, \mathbf z\right\|_2+
	2\sum_{k=1}^2\left\|\nabla \mathbf z_k\right\|_2
	\left\|\nabla \mathbf y_k\right\|_\infty\left\|\mathbf{curl}\,
	\mathbf z\right\|_2\right) \vspace{2mm}\\
	&\leq \left( S_4^2\left|y\right|_{H^1}+\kappa\alpha 
	\left\|y\right\|_{H^3}\right) \|\nabla z\|_2^2
	\end{array}$$
and the claimed result is proved.$\hfill\Box$
%%%%%%%%%%%%%%%%%%%%%%%%%%%%%%%%%%%%%%%%%%%%%%%%%%%%%%%
%%%%%%%%%%%%%%%%%%%%%%%%%%%%%%%%%%%%%%%%%%%%%%%%%%%%%%%
\section{State equation}
\setcounter{equation}{0} \label{state_equation_section}
\subsection{Existence and uniqueness results for the state equation}
\label{existence_section}
%%%%%%%%%%%%%%%%%%%%%%%%%%%%%%%%%%%%%%%%%%%%%%%%%%%%%%%
%%%%%%%%%%%%%%%%%%%%%%%%%%%%%%%%%%%%%%%%%%%%%%%%%%%%%%%
  The state equation can be written in a variational form by taking its scalar product with a test function in $ V$.
%%%%%%%%%%%%%%%%%%%%%%%%%%%%%%%%%%%%%%%%%%%%%%%%%%%%%%%
\begin{definition} Let $ u\in L^2(\Omega)$. A function $ y\in  V_2$ is a solution of $(\ref{equation_etat})$ if
	\begin{equation}\label{var_form_state}\nu\left(\nabla y,\nabla \phi\right)+\left( \mathbf{curl}\, \sigma(\mathbf y)
	\times  \mathbf y, \boldsymbol\phi\right)=\left( u, \phi\right) \qquad \mbox{for all} \  \phi\in  V.\end{equation}
\end{definition}
Due to Lemma \ref{non_lin_curl}, the nonlinear term in the previous definition can be understood in the following sense
	$$\begin{array}{ll}\left( \mathbf{curl}\, \sigma(\mathbf y)\times  \mathbf y, \boldsymbol\phi\right)&=
	b\left( \phi, y, \sigma(y)\right)-b\left( y,
	 \phi, \sigma(y)\right)
	\vspace{2mm}\\&=b\left( y, y, \phi\right)
	-\alpha\left(b\left( \phi, y,\Delta  y\right)-b\left( y,
	 \phi,\Delta  y\right)\right).\end{array}$$
%%%%%%%%%%%%%%%%%%%%%%%%%%%%%%%%%%%%%%%%%%%%%%%%%%%%%%%
Equation (\ref{equation_etat}) was first studied by Cioranescu and Ouazar (\cite{O81}, \cite{CO84}) in the case of Dirichlet boundary conditions and simply connected domains. These authors proved existence and uniqueness of solutions by using Galerkin's method in the basis of the eigenfunctions of the operator $\mathrm{curl}\left(\mathrm{curl}\, \sigma(y)\right)$. More precisely, by using the fact that the imbedding of $V_2\subset V$ is compact, they prove the existence of a sequence of eigenfunctions $\left(e_j\right)_j\subset V_2$ corresponding to a sequence of eigenvalues $\left(\lambda_j\right)_j$ such that
	\begin{equation}\label{eigenfunc}
	\left(e_j,\phi\right)_{V_2}=
	\lambda_j\left(
	\left( e_j,\phi\right)
	+\alpha\left(\nabla e_j,\nabla\phi\right)\right),
	\qquad \mbox{for all} \ \phi\in V_2\end{equation}
with 
	$$0<\lambda_1<\cdots<\lambda_k<\cdots \longrightarrow +\infty.$$
The functions $ e_j$ form an orthonormal basis in $V$ and an orthogonal basis in $V_2$. Moreover, 
	$$ e_j\in H^4(\Omega), \quad 
	{\rm curl} \, \sigma(e_j)\in H^1(\Omega)$$
and 
	\begin{equation}\label{g_vs_curlg}\left({\rm curl}\, g,
	{\rm curl}\, \sigma(e_j)\right)=\lambda_j\left(g, e_j\right)\qquad \mbox{for all} \ g\in H({\rm curl};\Omega).\end{equation}
This method, designed to decompose the problem into a Stokes-like system for the velocity $y$ and a transport equation for $\mathrm{curl}\,\sigma(y)$, allows to establish the  existence of global solutions with  $H^3$ regularity in the  two dimensional case, and uniqueness and local existence in the three dimensional case. It has been  extented by Cioranescu and Girault \cite{CG97} to prove global existence in time in the three dimensional case and by Busuioc and Ratiu \cite{BR03} to study the case of Navier-slip boundary conditions. \vspace{2mm}\\
The following result deals with existence of a solution and is well known (see {\it e.g.} \cite{GS99}). For the convenience of the reader, the corresponding estimates are derived herafter. 
%%%%%%%%%%%%%%%%%%%%%%%%%%%%%%%%%%%%%%%%%%%%%%%%%%%%%%%
\begin{proposition} \label{existence_state}Let $ u\in H( \mathrm{curl};\Omega)$. Then problem $(\ref{equation_etat})$ admits at least one solution $ y\in  V_2$ and this solution satisfies the following estimates
	\begin{equation}\label{state_est1}
	\left|y\right|_{H^1}\leq \tfrac{S_2}{\nu}
	\left\| u\right\|_2,\end{equation}
	\begin{equation}\label{state_est2}
	\left|y\right|_{V_2}\leq \tfrac{1}{\nu}
	\left(S_2\left\| u\right\|_2+\alpha
	\left\| \mathrm{curl} \, u\right\|_2\right),\end{equation}
	\begin{equation}\label{state_est3}
	\left\| y\right\|_{H^3}\leq 
	\tfrac{\kappa}{\alpha\nu}\left(S_2\left\| u\right\|_2+\alpha
	\left\| \mathrm{curl} \, u\right\|_2\right),
	\end{equation}
where $\kappa$ is a positive constant depending only on $\Omega$.
\end{proposition}
%%%%%%%%%%%%%%%%%%%%%%%%%%%%%%%%%%%%%%%%%%%%%%%%%%%%%%%
{\bf Proof.}  Setting $\phi=y$ in (\ref{var_form_state}) and using the Poincar\'e inequality, we obtain
	$$\begin{array}{ll}\nu\left|y\right|_{H^1}^2&=
	\left(u,y\right)-
	\left( \mathbf{curl}\, \sigma\mathbf (\mathbf y)\times  \mathbf y, \mathbf y\right)
	=\left(u,y\right)\vspace{2mm}\\
	&\leq
	 \|u\|_2\|y\|_2\leq S_{2}\|u\|_2\left|y\right|_{H^1}.
	\end{array}$$
which gives (\ref{state_est1}). On the other hand, by applying the curl to 
$(\ref{equation_etat})$, we obtain
	$$-\nu \Delta \left(\mathrm{curl} \, y\right)+y\cdot \nabla
	 \mathrm{curl} \,\sigma(y)=\mathrm{curl} \, u$$
yielding
	\begin{equation}\label{transport_state}
	\mathrm{curl} \, \sigma(y)+\tfrac{\alpha}{\nu}\,
	y\cdot \nabla\left(\mathrm{curl} \,\sigma(y)\right)
	=\tfrac{\alpha}{\nu}\,
	\mathrm{curl} \, u+\mathrm{curl} \, y.\end{equation}
Multiplying by $\mathrm{curl} \, \sigma(y)$ and integrating, we get
	$$\begin{array}{ll}\left|y\right|_{V_2}^2&=\left\|\mathrm{curl}\, \sigma(y)\right\|_{2}^2\vspace{2mm}\\
	&=-\tfrac{\alpha}{\nu}\left(y\cdot \nabla\left(\mathrm{curl} \,\sigma(y)\right),\mathrm{curl}\, \sigma(y)\right)+\left(\tfrac{\alpha}{\nu}\,\mathrm{curl} \, u+\mathrm{curl} \, y,\mathrm{curl}\, 
	\sigma(y)\right)\vspace{2mm}\\
	&=\left(\tfrac{\alpha}{\nu}\,\mathrm{curl} \, u+\mathrm{curl} \, y,\mathrm{curl}\,
	 \sigma(y)\right)	\vspace{2mm}\\
	&\leq\left(\tfrac{\alpha}{\nu}\left\|\mathrm{curl} \, u\right\|_2+\left\|\mathrm{curl} \, y\right\|_2\right)
	\left|y\right|_{V_2}\vspace{2mm}\\
	&\leq\left(\tfrac{\alpha}{\nu}\left\|\mathrm{curl} \, u\right\|_2+\left|y\right|_{H^1}\right)
	\left|y\right|_{V_2}
	\end{array}$$
and thus
	$$\left|y\right|_{V_2}\leq \left|y\right|_{H^1}
	+\tfrac{\alpha}{\nu}
	\left\| \mathrm{curl} \, u\right\|_2.$$
This estimate together with (\ref{state_est1}) gives (\ref{state_est2}).
Finally, since $\mathrm{curl}\,\Delta y\in L^2(\Omega)$ and $\nabla\cdot \left(\mathrm{curl}\,\Delta y\right)=0$, by arguing as in the proof of Lemma \ref{infty_V2}, we can establish the existence of a unique function $\psi\in H^1(\Omega)$ such that 
	$$-\Delta y -\psi+\nabla \pi=0$$
and 
	\begin{equation}\label{y2_sigma} \left\|y\right\|_{H^3}
	\leq c\left\|\psi\right\|_{H^1}
	\leq \kappa\left\|\mathrm{curl}\,\Delta y\right\|_{2},\end{equation}
where $\kappa$ is a positive constant only depending on $\Omega$.
Combining (\ref{transport_state}) and $(\ref{y2_sigma})$, we obtain
	$$\begin{array}{ll}\left\|y\right\|^2_{H^3}&\leq 
	\kappa^2\left\|\mathrm{curl} \, \Delta y\right\|_2^2\vspace{2mm}\\
	&=\left(\tfrac{\kappa}{\alpha}\right)^2\left\|\mathrm{curl} \, 
	y-\mathrm{curl} \, \sigma(y)\right\|_2^2\vspace{2mm}\\
	&=\left(\tfrac{\kappa}{\alpha}\right)^2\left(
	\left\|\mathrm{curl} \, y\right\|_2^2
	+\left|y\right|_{V_2}^2-2
	\left(\mathrm{curl} \, y,\mathrm{curl} \, \sigma(y)\right)
	\right)\vspace{2mm}\\
	&=\left(\tfrac{\kappa}{\alpha}\right)^2\left(
	\left\|\mathrm{curl} \, y\right\|_2^2
	-\left|y\right|_{V_2}^2
	+\tfrac{2\alpha}{\nu}
	\left(\mathrm{curl} \, u,\mathrm{curl} \, \sigma(y)\right)
	\right)\vspace{2mm}\\
	&\leq \left(\tfrac{\kappa}{\alpha}\right)^2\left(
	\left\|\mathrm{curl} \, y\right\|_2^2
	+\left(\tfrac{\alpha}{\nu}
	\left\|\mathrm{curl} \, u\right\|_2
	\right)^2\right)\vspace{2mm}\\
	&\leq \left(\tfrac{\kappa}{\alpha}\right)^2\left(
	\left|y\right|_{H^1}^2
	+\left(\tfrac{\alpha}{\nu}
	\left\|\mathrm{curl} \, u\right\|_2
	\right)^2\right)\end{array}$$
and thus
	$$\left\|y\right\|_{H^3}\leq \tfrac{\kappa}{\alpha}
	\left(\left|y\right|_{H^1}+\tfrac{\alpha}{\nu}
	\left\|\mathrm{curl} \, u\right\|_2\right).$$
Estimate (\ref{state_est3}) is then a direct consequence of (\ref{state_est1}).$\hfill\Box$\vspace{2mm}\\
%%%%%%%%%%%%%%%%%%%%%%%%%%%%%%%%%%%%%%%%%%%%%%%%%%%%%%%
As in the case of Navier-Stokes equations, uniqueness of the solution is guaranteed under a restriction on the data. Additional regularity of the solution is obtained under the same restriction for more regular data.  
%%%%%%%%%%%%%%%%%%%%%%%%%%%%%%%%%%%%%%%%%%%%%%%%%%%%%%%
\begin{proposition}\label{uniqueness_state} 
Assume that $ u\in H(\mathrm{curl};\Omega)$. There exists a positive constant $ \bar\kappa$, depending only on $\Omega$, such that if $ u$ satisfies 
	\begin{equation}\label{uniqueness_condition}
	 \bar\kappa\left(\left\| u\right\|_2+\alpha
	\left\| \mathrm{curl} \, u\right\|_2\right)<\nu^2,\end{equation}
then equation $(\ref{equation_etat})$ admits a unique solution $y$. Moreover, if ${\rm curl} \, u \in H^1(\Omega)$ then $y\in H^4(\Omega)$ and the following estimate holds
$$\left(1-\tfrac{\bar\kappa}{\nu^2}\left(
	\left\|u\right\|_2+
	\alpha\left\|{\rm curl} \,u\right\|_2\right)\right)
	\left|{\rm curl}\, \sigma(y) \right|_{H^1}
	\leq \tfrac{\kappa}{\alpha\nu}
	\left(\left\|u\right\|_2+
	\alpha\left\|{\rm curl} \,u\right\|_2+\alpha^2\left|
	{\rm curl} \, u \right|_{H^1}\right).$$
\end{proposition}
%%%%%%%%%%%%%%%%%%%%%%%%%%%%%%%%%%%%%%%%%%%%%%%%%%%%%%%
{\bf Proof.} Assume that $ y_1$ and $ y_2$ are two solutions of (\ref{equation_etat}) corresponding to $ u$ and denote by $ y$ the difference $ y_1- y_2$. By setting $\phi=y$ in the variational formulation (\ref{var_form_state}),  we deduce that
	$$\nu\left|y\right|_{H^1}^2+\left( \mathbf{curl} \,\sigma(\mathbf y_1)
	\times  \mathbf y_1- \mathbf{curl}\, \sigma(\mathbf y_2)
	\times  \mathbf y_2,\mathbf y\right)=0.$$
Observing that
	$$\mathbf{curl} \,\sigma(\mathbf y_1)
	\times  \mathbf y_1- \mathbf{curl}\, \sigma(\mathbf y_2)
	\times  \mathbf y_2=\mathbf{curl} \,\sigma(\mathbf y_1)
	\times  \mathbf y+\mathbf{curl} \,\sigma(\mathbf y)	\times \mathbf  y_2,$$
and taking into account  Lemma \ref{non_lin_curl}, we deduce that
	$$\nu\left|y\right|_{H^1}^2+\left( \mathbf{curl} \,\sigma(\mathbf y)
	\times  \mathbf y_2,\mathbf y\right)=0.$$
Due to  Lemma \ref{rm2}, (\ref{state_est1}) and
 (\ref{state_est2}), it follows that
	$$\begin{array}{ll}\left|y\right|_{H^1}^2
	&\leq \tfrac{1}{\nu}\left( S_4^2\left|y_2\right|_{H^1}+
	\kappa\alpha\left\| y_2
	\right\|_{H^3}\right)\left|y\right|_{H^1}^2\vspace{2mm}\\
	&\leq \tfrac{1}{\nu^2}\left( S_4^2 S_2\left\|u\right\|_2+
	\kappa\left(\left\| u\right\|_2+\alpha
	\left\| \mathrm{curl} \, u\right\|_2\right)\right)
	\left|y\right|_{H^1}^2\vspace{2mm}\\
	&\leq \tfrac{ \bar\kappa}{\nu^2}\left(\left\| u\right\|_2+\alpha
	\left\| \mathrm{curl} \, u\right\|_2\right)
	\left|y\right|_{H^1}^2
	\end{array}$$
implying that $ y_1= y_2$ if condition (\ref{uniqueness_condition}) is satisfied. This proves the uniqueness result. The regularity result can be similarly established by using classical arguments on the transport equation. For the convenience of the reader, we will give a sketch of the proof and only derive the estimate that shall be applied to the solution of a Galerkin approximation of the problem.
By taking the gradient in (\ref{transport_state}), we can see that $\varphi =\nabla\left(\mathrm{curl}\, \sigma(y)\right) $ is the solution of the following transport equation
	$$\varphi +\tfrac{\alpha}{\nu}\,
	 y  \cdot \nabla \varphi +\tfrac{\alpha}{\nu}
	\left(\nabla y\right)^\top \cdot \varphi=
	\nabla\left(\tfrac{\alpha}{\nu}\,{\rm curl} \, 
	u +{\rm curl} \, y \right).$$
Therefore
	$$\begin{array}{ll}\left\|\varphi \right\|_2^2&
	=\left(\nabla\left(\tfrac{\alpha}{\nu}\,{\rm curl} \, 
	u +{\rm curl} \, y \right),
	\varphi \right)
	-\left(\tfrac{\alpha}{\nu}\left(\nabla y 
	\right)^\top \cdot
	 \varphi,\varphi\right)\vspace{1mm}\\
	&\leq \left(\tfrac{\alpha}{\nu} \left\|{\rm curl} \, u
	\right\|_{H^1}+\left\|{\rm curl} \, y
	 \right\|_{H^1}\right)
	\left\|\varphi \right\|_2+
	\tfrac{\alpha}{\nu}\left\|\nabla y \right\|_\infty
	 \left\|\varphi \right\|_2^2\vspace{2mm}\\
	&\leq \left(\tfrac{\alpha}{\nu} \left\|{\rm curl} \, u
	\right\|_{H^1}+c\left\|y
	 \right\|_{H^2}\right)
	\left\|\varphi \right\|_2+
	\tfrac{c\alpha}{\nu}\left\|y \right\|_{H^3}
	 \left\|\varphi \right\|_2^2\vspace{2mm}\\
	&\leq \left(\tfrac{\alpha}{\nu} \left\|{\rm curl} \, u
	\right\|_{H^1}+\tfrac{c}{\alpha}
	\left\|{\rm curl} \, \sigma(y) \right\|_{2}\right)
	\left\|\varphi \right\|_2+
	\tfrac{ \bar\kappa}{\nu^2}\left(\left\| u\right\|_2+\alpha
	\left\| \mathrm{curl} \, u\right\|_2\right)
	\left\|\varphi \right\|_2^2\vspace{2mm}\\
	&\leq \tfrac{\kappa}{\alpha\nu}
	\left(\left\|u\right\|_2+
	\alpha\left\|{\rm curl} \,u\right\|_2+\alpha^2\left\|\nabla\left({\rm curl} \, u \right)\right\|_2\right)
	\left\|\varphi \right\|_2+
	\tfrac{ \bar\kappa}{\nu^2}\left(\left\| u\right\|_2+\alpha
	\left\| \mathrm{curl} \, u\right\|_2\right)
	\left\|\varphi \right\|_2^2.
\end{array}$$
This gives the estimate and shows that $\mathrm{curl}\, \sigma(y)$ belongs to $H^1(\Omega)$. Arguing as in the proof of (\ref{yh_curl_sigma}) and (\ref{y2_sigma}), it follows that $y\in H^4(\Omega)$.  $\hfill\Box$
%%%%%%%%%%%%%%%%%%%%%%%%%%%%%%%%%%%%%%%%%%%%%%%%%%%%%%%
\begin{remark} \label{remark1} Notice that
	$ \bar\kappa >S_4^2 S_2$. This  
implies that if $u$ satisfies the condition stated in the previous proposition, then the corresponding Navier-Stokes equation has a unique weak solution. 
\end{remark}

%%%%%%%%%%%%%%%%%%%%%%%%%%%%%%%%%%%%%%%%%%%%%%%%%%%%%%%
%%%%%%%%%%%%%%%%%%%%%%%%%%%%%%%%%%%%%%%%%%%%%%%%%%%%%%%
\subsection{Linearized state equation} \label{lin_section}
%\setcounter{equation}{0}
%%%%%%%%%%%%%%%%%%%%%%%%%%%%%%%%%%%%%%%%%%%%%%%%%%%%%%%
%%%%%%%%%%%%%%%%%%%%%%%%%%%%%%%%%%%%%%%%%%%%%%%%%%%%%%%
The aim of this section is to study the solvability, in an adequate setting, of the linearized equation associated to the nonlinear state equation. Its solution is involved in the definition of the directional derivative of the control-to-state mapping and is related, through a suitable Green formula, to the adjoint state.\vspace{2mm}\\
Let  $u\in H( \mathrm{curl};\Omega)$, let $ y\in V_2$ be a corresponding solution of $(\ref{equation_etat})$ and consider the linear equation 
\begin{equation}\label{linearized}\left\{
  \begin{array}{ll}
    -\nu\Delta  \mathbf z+ \mathbf{curl}\,\sigma(\mathbf z)\times \mathbf y+
	 \mathbf{curl}\,\sigma(\mathbf y)\times \mathbf z+\nabla \pi=\mathbf w&\quad\mbox{in} \ \Omega,\vspace{2mm}\\
	 \mathrm{div} \,  \mathbf z=0&\quad\mbox{in} \ \Omega,\vspace{2mm}\\
    \mathbf z=0&\quad\mbox{on}\ \Gamma,
  \end{array}
\right .\end{equation}
where $w\in L^2(\Omega)$.
%%%%%%%%%%%%%%%%%%%%%%%%%%%%%%%%%%%%%%%%%%%%%%%%%%%%%%%	
\begin{definition}  A function $ z\in  V_2$ is a solution of $(\ref{linearized})$ if
	\begin{equation}\label{var_lin}
	\nu\left(\nabla z,\nabla \phi\right)+
	\left(\mathbf{curl}\, \sigma(\mathbf z)\times
	  \mathbf y+\mathbf{curl}\, \sigma(\mathbf y)\times  \mathbf z, 
	\boldsymbol\phi\right)=
	\left( w, \phi\right) \qquad \mbox{for all} \ 
	\phi\in  V.\end{equation}
\end{definition}
%%%%%%%%%%%%%%%%%%%%%%%%%%%%%%%%%%%%%%%%%%%%%%%%%%%%%%%%
In analogy to the state equation, by taking into account Lemma 
\ref{non_lin_curl}, we can rewrite the previous variational formulation as follows:
	$$\nu\left(\nabla z,\nabla \phi\right)+
	b\left( \phi, y, \sigma(z) \right)-b\left( y, \phi, \sigma(z) \right)+
	b\left( \phi, z, \sigma(y) \right)-b\left( z, \phi, \sigma(y) \right)=
	\left( w, \phi\right)$$
for all $\phi\in V$.\vspace{2mm}\\
As already mentioned, the special Galerkin basis used to study the state equation (\ref{equation_etat}) is particularly well adapted and allows to prove existence of regular solutions with minimal assumptions on the data. 
Seeming appropriate, the application of the same arguments to study the solvability of the linearized equation (\ref{linearized}) leads, however, to additional, yet expectable, issues. Indeed, after deriving the $H^1$ a priori estimate, this technique will naturally {\it imposes} the derivation of a $L^2$ estimate for 
$\mathrm{curl}\, \sigma(z)$  (and thus $H^3$ for $z$). This term should satisfy the transport equation
	$$\mathrm{curl}\, \sigma(z)+\tfrac{\alpha}{\nu} \,y\cdot \nabla 
	\left(\mathrm{curl}\, \sigma(z)\right)
	+\tfrac{\alpha}{\nu}\,z\cdot \nabla 
	\left(\mathrm{curl}\, \sigma(y)\right)=\tfrac{\alpha}{\nu}
	\,\mathrm{curl}\, w+\mathrm{curl}\, z$$
and in order to obtain the desired estimate,  we need to guarantee that the coefficient $\mathrm{curl}\,\sigma(y)$ appearing in the linearized operator belongs to $H^1(\Omega)$. Following Proposition \ref{uniqueness_state}, this can be achieved if we consider more regular data in the state equation and impose additional restrictions on their size.\vspace{1mm}\\
On the other hand, let us observe that the variational formulation stated above is well defined for $\sigma(z)\in L^2(\Omega)$ (and thus for $z\in H^2(\Omega)$) and that this regularity would be sufficient to carry out our analysis and derive the necessary optimality conditions. We might consider less restrictive choices for the Galerkin basis, but technical difficulties inherent to Dirichlet boundary conditions need to be managed. 
Formally, the natural way to obtain the $H^2$ a priori estimates would be to multiply (\ref{linearized}) by $\sigma(z)$
and to integrate. The main difficulty is then to deal with the pressure term
	$$\left(\nabla \pi,\sigma(z)\right)=-\left(\pi,\mathrm{div}\, \sigma(z)\right)+\int_\Gamma \pi n\cdot \sigma(z)=\int_\Gamma \pi n\cdot \sigma(z)$$
 that does not vanish, unless  $\sigma(z)$ is tangent to the boundary, and that we do no know how to adequately estimate. \vspace{2mm}\\
The next result deals with existence of a regular solutions of the linearized equation when the corresponding data are sufficiently regular.
%%%%%%%%%%%%%%%%%%%%%%%%%%%%%%%%%%%%%%%%%%%%%%%%%%%%%%%
\begin{proposition} \label{ex_uniq_lin}Let $u\in H( \mathrm{curl};\Omega)$ satisfying condition $(\ref{uniqueness_condition})$ and such that  ${\rm curl} \, u \in H^1(\Omega)$,  and let $y\in V_2\cap H^4(\Omega)$ be the corresponding solution of $(\ref{equation_etat})$. Then equation $(\ref{linearized})$ admits a unique solution $ z\in V_2$. Moreover, the following estimates hold
	\begin{equation}\label{est_H1_z}	
	\left(1-\tfrac{\bar\kappa}{\nu^2}\left(\left\| u\right\|_2+\alpha
	\left\| \mathrm{curl} \, u\right\|_2\right)\right)\left|z\right|_{H^1}
	\leq \tfrac{S_2}{\nu} \left\| w\right\|_2,\end{equation}
	\begin{equation}\label{est_V2_z}\left|z\right|_{V_2}
	\leq \tfrac{2\alpha}{\nu}\left\|{\rm curl}\, w\right\|_2+
	 \kappa\left(\tfrac{\alpha}{\nu^\frac{3}{2}}
	\left|\mathrm{curl}\,\sigma(y)\right|_{H^1}^\frac{3}{2}+
	1\right)\left|z\right|_{H^1},\end{equation} 
where $\kappa$ is a positive constant depending only on $\Omega$.
\end{proposition}
%%%%%%%%%%%%%%%%%%%%%%%%%%%%%%%%%%%%%%%%%%%%%%%%%%%%%%%
%%%%%%%%%%%%%%%%%%%%%%%%%%%%%%%%%%%%%%%%%%%%%%%%%%%%%%%
%\subsubsection{Galerkin method and estimates for approximate solutions}
%%%%%%%%%%%%%%%%%%%%%%%%%%%%%%%%%%%%%%%%%%%%%%%%%%%%%%%
%%%%%%%%%%%%%%%%%%%%%%%%%%%%%%%%%%%%%%%%%%%%%%%%%%%%%%%
{\bf Proof.} The proof of Proposition \ref{ex_uniq_lin} is split into three steps. We first establish the existence of an approximate solution and a first estimate in $ H^1(\Omega)$. Next, we derive an estimate in $H^3(\Omega)$ and then we pass to the limit. \vspace{2mm}\\
The solution of $(\ref{linearized})$ is constructed by means of Galerkin's discretization, by expanding the linearized state $z$ in the basis introduced in the previous section.
The approximate problem is defined by
	\begin{equation}\label{faedo_galerkin_lin}
	\left\{\begin{array}{ll}\mbox{Find} \ 
	 z_m=\displaystyle\sum_{j=1}^m \zeta_{j}
	  e_j \ \mbox{solution}, \ \mbox{for} \ 1\leq j\leq m, \ 
	\mbox{of}\vspace{3mm}\\
	\nu\left(\nabla z_m,\nabla e_j\right)+
	\left( \mathbf{curl}\,\sigma(\mathbf z_m)\times  \mathbf y+
	 \mathbf{curl}\,\sigma(\mathbf y)\times
	  \mathbf z_m, \mathbf  e_j\right)=
	\left(w,  e_j\right).\end{array}\right.
	\end{equation}
%%%%%%%%%%%%%%%%%%%%%%%%%%%%%%%%%%%%%%%%%%%%%%%%%%%%%%%
{\it Step 1. Existence of the discretized solution and a priori $H^1$ estimate.} We prove that the $H^1$ estimate can be derived if $u$ satisfies the condition $(\ref{uniqueness_condition})$. Let $m$ be fixed and consider $P: \ \mathbb{R}^m \longrightarrow \mathbb{R}^m$ defined by
	$$\left(P \zeta\right)_j=\nu\left(\nabla z_m,\nabla e_j\right)+\left(\mathbf{curl}\,\sigma(\mathbf z_m)\times \mathbf y+
	\mathbf{curl}\, \sigma(\mathbf y)\times \mathbf z_m,  \mathbf  e_j\right)-
	\left( w,   e_j\right),	$$
where $ z_m=\sum_{j=1}^m\zeta_j    e_j$. The mapping $P$ is obviously continuous. Let us prove that $P(\zeta)\cdot \zeta>0$ if $|\zeta|$ is sufficiently large. Classical arguments together with Lemma \ref{rm2} yields
	\begin{align}\label{P_zeta}{P}(\zeta)\cdot \zeta&=
	\nu \left|z_m\right|_{H^1}^2+\left(\mathbf{curl}\,\sigma(\mathbf z_m)\times  \mathbf y, \mathbf  z_m\right)-
	\left( w, z_m\right)\\
	&\geq \left(\nu-\left( S_4^2\left|y\right|_{H^1}+
	\kappa\alpha\| y\|_{H^3}\right)\right)
	\left|z_m\right|_{H^1}^2
	-\left\| w\right\|_2\left\|z_m\right\|_2\nonumber\\
	&\geq \left(\nu-\left( S_4^2\left|y\right|_{H^1}+
	\kappa\alpha\| y\|_{H^3}\right)\right)
	 \left|\zeta\right|^2-\left\| w\right\|_2
	\left|\zeta\right|\nonumber\\
	&\geq \left(\nu-\tfrac{\bar\kappa}{\nu}
	\left(\left\| u\right\|_2+\alpha
	\left\| \mathrm{curl} \, u\right\|_2\right)
	\right)\left|\zeta\right|^2-
	\left\| w\right\|_2\left|\zeta\right|\nonumber\\
	&\longrightarrow +\infty \qquad \mbox{when} \ 
	|\zeta|\rightarrow +\infty.\nonumber\end{align}
Due to the Brouwer theorem, we deduce that there exists $ \zeta^\ast\in \mathbb{R}^m$ such that ${P}\left( \zeta^\ast\right)=0$ and thus $ z_m=\sum_{j=1}^m\zeta_j^\ast    e_j$ is a solution of problem (\ref{faedo_galerkin_lin}). Due to $(\ref{P_zeta})$ and  Lemma \ref{rm2}, it follows that
	$$\begin{array}{ll}\nu \left|z_m\right|_{H^1}^2&=\left( \mathbf{curl}\,\sigma(\mathbf z_m)\times  \mathbf y,  \mathbf z_m\right)-
	\left( w, z_m\right)\vspace{2mm}\\
	&\leq \left( S_4^2\left|y\right|_{H^1}+
	\kappa\alpha
	\left\| y\right\|_{H^3}\right)
	\left|z_m\right|_{H^1}^2+
	S_2\left\| w\right\|_2
	 \left|z_m\right|_{H^1}\vspace{2mm}\\
	&\leq\tfrac{\bar\kappa}{\nu}\left(\left\| u\right\|_2+\alpha
	\left\| \mathrm{curl} \, u\right\|_2\right)
	\left|z_m\right|_{H^1}^2+
	S_2\left\| w\right\|_2
	 \left|z_m\right|_{H^1}
	\end{array}$$
which gives 
\begin{equation}\label{est_H1_zm}\left(
	1-\tfrac{ \bar\kappa}{\nu^2}
	\left(\left\| u\right\|_2+\alpha
	\left\| \mathrm{curl} \, u\right\|_2\right)\right)\left|z_m\right|_{H^1}
	\leq \tfrac{S_2}{\nu} \left\| w\right\|_2.\end{equation}
%%%%%%%%%%%%%%%%%%%%%%%%%%%%%%%%%%%%%%%%%%%%%%
{\it Step 2. A priori $H^3$ estimate.} By taking into account (\ref{eigenfunc}) and (\ref{g_vs_curlg}) we have
	$$\begin{array}{ll}\left|z_m\right|_{V_2}^2&=\left\|\mathrm{curl} \,\sigma(z_m)\right\|_2^2\vspace{2mm}\\
	&=\displaystyle\sum_{j=1}^m \zeta_j^\ast\left(\mathrm{curl} \,\sigma(z_m),{\rm curl}\, \sigma(e_j)\right)\vspace{1mm}\\
	&=\displaystyle\sum_{j=1}^m \zeta_j^\ast\left(\lambda_j-1\right)
	\left(\left(z_m, e_ j \right)+\alpha
	\left(\nabla z_m,\nabla e_ j \right)\right)
	\vspace{1mm}\\
	&=\displaystyle\sum_{j=1}^m \zeta_j^\ast\left(\lambda_j-1\right)
	\left(\mathbf z_m-\tfrac{\alpha}{\nu}\left(\mathbf{curl} \,\sigma(\mathbf z_m)\times \mathbf y+
	\mathbf{curl}\, \sigma\left(\mathbf y\right)\times \mathbf z_m-\mathbf w\right),
	 \mathbf e_ j \right)\vspace{1mm}\\
	&=\displaystyle\sum_{j=1}^m \zeta_j^\ast
	\left({\rm curl}\,z_m+\tfrac{\alpha}{\nu}\, {\rm curl}\,w,
	{\rm curl}\, \sigma(e_ j)\right)\vspace{1mm}\\
	&-\tfrac{\alpha}{\nu}\displaystyle\sum_{j=1}^m \zeta_j^\ast
	\left(\mathbf{curl}\left(\mathbf{curl} \,\sigma(\mathbf z_m)\times \mathbf y+
	\mathbf{curl}\, \sigma\left(\mathbf y\right)\times \mathbf z_m\right),
	\mathbf{ curl}\, \sigma(\mathbf e_j)\right)
	\end{array}$$
yielding
$$\begin{array}{ll}\left|z_m\right|_{V_2}^2
	&=\left({\rm curl}\,z_m+\tfrac{\alpha}{\nu}\,
	 {\rm curl}\,w,
	{\rm curl}\, \sigma(z_m)\right)\vspace{3mm}\\
	&-\tfrac{\alpha}{\nu}
	\left(\mathbf{curl}\left(\mathbf{curl} \,\sigma(\mathbf z_m)\times \mathbf y+
	\mathbf{curl}\, \sigma\left(\mathbf y\right)\times \mathbf z_m\right),
	\mathbf{ curl}\, \sigma(\mathbf z_m)\right)\vspace{2mm}\\
	&=
	\left({\rm curl}\,z_m+\tfrac{\alpha}{\nu}\,
	 {\rm curl}\,w,	{\rm curl}\, \sigma(z_m)\right)\vspace{2mm}\\
	&-\tfrac{\alpha}{\nu}\left(b\left(y,{\rm curl}\, \sigma(z_m),
	{\rm curl}\, \sigma(z_m)\right)
	+b\left(z_m,{\rm curl}\, \sigma(y),{\rm curl}\,
	 \sigma(z_m)\right)\right)\vspace{2mm}\\
	&=\left({\rm curl}\,z_m+\tfrac{\alpha}{\nu}\,
	 {\rm curl}\,w,	{\rm curl}\, \sigma(z_m)\right)
	-\tfrac{\alpha}{\nu}\,b\left(z_m,{\rm curl}\, 
	\sigma(y),{\rm curl}\, \sigma(z_m)\right).
	\end{array}$$
Due to Lemma \ref{infty_V2}, we deduce that
$$\begin{array}{ll}\left|z_m\right|_{V_2}&
	\leq\left\|{\rm curl}\,z_m\right\|_2+
	\tfrac{\alpha}{\nu}
	\left\|{\rm curl}\, w\right\|_2+
	\tfrac{\alpha}{\nu}\left\|z_m\right\|_\infty
	\left|\mathrm{curl}\,\sigma(y)\right|_{H^1}
	\vspace{2mm}\\
	&\leq 
	\left\|{\rm curl}\,z_m\right\|_2+
	\tfrac{\alpha}{\nu}
	\left\|{\rm curl}\, w\right\|_2+
	\tfrac{c\alpha^\frac{2}{3}}{\nu}\left|z_m\right|_{H^1}^{\frac{2}{3}}\left|z_m\right|_{V_2}^{\frac{1}{3}}
	\left|{\rm curl}\, \sigma(y)\right|_{H^1}\end{array}$$
and by using the Young inequality we finally obtain
	\begin{equation}\label{est_V2_zm}
	\left|z_m\right|_{V_2}
	\leq \tfrac{2\alpha}{\nu}\left\|{\rm curl}\, w\right\|_2+
	 \kappa\left(\tfrac{\alpha}{\nu^\frac{3}{2}}
	\left|\mathrm{curl}\,\sigma(y)\right|_{H^1}^\frac{3}{2}+
	1\right)\left|z_m\right|_{H^1}.\end{equation}
{\it Step 3. Passing to the limit.} It remains to pass to the limit with respect to $m$. From estimates (\ref{est_H1_zm}) and (\ref{est_V2_zm}), it follows that if $u$ satisfies condition $(\ref{uniqueness_condition})$ then there exists a subsequence, still indexed by $m$, and function $ z\in  V_2$ such that
	$$ z_m \longrightarrow  z \qquad \mbox{weakly in} \  V_2.$$
By passing to the limit in (\ref{faedo_galerkin_lin}), we obtain for every $j\geq 1$
	$$\nu\left(\nabla z,\nabla e_j\right)+
	b\left(   e_j, y, \sigma(z)\right)-b\left( y,   e_j,
	\sigma(z)\right)+b\left(   e_j, z, \sigma(y) \right)
	-b\left( z,   e_j, \sigma(y)\right)=
	\left( w,   e_j\right)	$$
and by density we prove that $ z$ satisfies the variational formulation. Moreover, $z$ satisfies estimates (\ref{est_H1_z}) and (\ref{est_V2_z}). Finally, since (\ref{linearized}) is linear, the uniquess result is direct consequence of estimate (\ref{est_H1_z}).$\hfill\Box$
%%%%%%%%%%%%%%%%%%%%%%%%%%%%%%%%%%%%%%%%%%%%%%%%%%%%%%%
%%%%%%%%%%%%%%%%%%%%%%%%%%%%%%%%%%%%%%%%%%%%%%%%%%%%%%%
\section{Analysis of the control-to-state mapping} \label{sec_lip}
\setcounter{equation}{0}
\subsection{Sequential and Lipschitz continuity}
%%%%%%%%%%%%%%%%%%%%%%%%%%%%%%%%%%%%%%%%%%%%%%%%%%%%%%%
%%%%%%%%%%%%%%%%%%%%%%%%%%%%%%%%%%%%%%%%%%%%%%%%%%%%%%%
We are first concerned with continuity properties of the control-to-state mapping in adequate topologies. 
%%%%%%%%%%%%%%%%%%%%%%%%%%%%%%%%%%%%%%%%%%%%%%%%%%%%%%%
\begin{proposition} \label{continuity_H3} Let $U$ be a bounded closed subset of $H(\mathrm{curl};\Omega)$. Then the control-to-state mapping is sequentially continuous from $U$, endowed with its weak topology, into 
$H^2(\Omega)$.
%Let $(u_k)_k$ be a bounded sequence in $H(\mathrm{curl};\Omega)$ and let $y_k$ be a solution of $(\ref{equation_etat})$ corresponding to $u_k$. There exists a subsequence, still indexed by $k$, such that $(u_k)_k$ weakly converges to some $u$ in $H(\mathrm{curl};\Omega)$ and $(y_k)_k$ strongly converges in $H^3(\Omega)$ to a solution of $(\ref{equation_etat})$ corresponding to $u$.
\end{proposition}
%%%%%%%%%%%%%%%%%%%%%%%%%%%%%%%%%%%%%%%%%%%%%%%%%%%%%%%
{\bf Proof.} Let $(u_k)_k \subset U$ be a sequence converging to $u$ in the weak topology of $H(\mathrm{curl};\Omega)$ and let $y_k$ be a solution of $(\ref{equation_etat})$ corresponding to $u_k$. 
Due to estimates (\ref{state_est2}) and (\ref{state_est3}), we have
	$$\left|y_k\right|_{V_2}\leq \tfrac{1}{\nu}
	\left(S_2\left\|u_k\right\|_2+\alpha 
	\left\|\mathrm{curl}\, u_k\right\|_2\right),$$
	$$\left\|y_k\right\|_{H^3}\leq \tfrac{\kappa}{\alpha \nu}
	\left(S_2\left\|u_k\right\|_2+\alpha 
	\left\|\mathrm{curl}\, u_k\right\|_2\right)$$
and since $(u_k)_k$ is uniformly bounded in $H(\mathrm{curl};\Omega)$, we deduce that the sequence $(y_k)_k$ is bounded in $V_2$. Then there exists a subsequence, still indexed by $k$, and $y\in V_2$, such that  $\left(y_k\right)_{k}$ weakly converges to $y$ in $H^3(\Omega)$ and (by using compactness results on Sobolev spaces) strongly in $H^2(\Omega)$. 
By passing to the limit in the variational formulation corresponding to $y_k$, we obtain
	$$\nu \left(\nabla y,\nabla \phi\right)+
	b\left( \phi, y, \sigma(y)\right)-b\left(y,
	 \phi, \sigma(y)\right)=(u,\phi) \qquad \mbox{for all} \ 
	\phi \in V$$
implying that $y$ is a solution of (\ref{equation_etat}) corresponding to $u$, and the claimed result is proved.$\hfill\Box$\vspace{2mm}\\
Next, we analyze the local Lipschitz continuity of the state with respect to the control variable. More precisely, if $u_1$, $u_2$ are two controls in $H({\rm curl};\Omega)$ and if $y_1$,  $y_2$ are two corresponding states then, by assuming that one of the control variables satisfies the restriction $(\ref{uniqueness_condition})$, we estimate $\left|y_1-y_2\right|_{H^1}$ with respect to $\left\|u_1-u_2\right\|_2$. Under the additional assumption that this control variable is regular enough, we can also estimate $\left|y_1-y_2\right|_{V_2}$ with respect to $\left\|u_1-u_2\right\|_{H({\rm curl};\Omega)}$.
%%%%%%%%%%%%%%%%%%%%%%%%%%%%%%%%%%%%%%%%%%%%%%%%%%%%%%%
\begin{proposition}
Let $ u_1,  u_2\in H(\mathrm{curl};\Omega)$ with $u_2$ satisfying condition $(\ref{uniqueness_condition})$, and let $ y_1,  y_2 \in  V_2$ be corresponding solutions of $(\ref{equation_etat})$. Then the following estimate holds
	\begin{equation}\label{lipschitz_est_H1}
	\left(1-\tfrac{ \bar\kappa}{\nu^2}\left(\left\| u_2\right\|_2+\alpha
	\left\| \mathrm{curl} \, u_2\right\|_2\right)\right)
	\left|y_1- y_2\right|_{H^1}\leq \tfrac{S_2}{\nu} \left\| u_1- u_2\right\|_2.
	\end{equation}
Moreover, if ${\rm curl}\,u_2$ belongs to $H^1(\Omega)$ then
	\begin{equation}\label{lipschitz_est_V2}
	\left|y_1-y_2\right|_{V_2}
	\leq \tfrac{2\alpha}{\nu}\left\|{\rm curl}\, (u_1-u_2)\right\|_2
	+\kappa\left(\tfrac{\alpha}{\nu^\frac{3}{2}}
	\left|{\rm curl}\, \sigma(y_2)\right|_{H^1}^\frac{3}{2}+
	1\right)\left|y_1-y_2\right|_{H^1},
	\end{equation}
where $\kappa$ is a positive constant only depending on $\Omega$.
\end{proposition}
%%%%%%%%%%%%%%%%%%%%%%%%%%%%%%%%%%%%%%%%%%%%%%%%%%%%%%%
{\bf Proof.} The proof is split into two steps.\vspace{2mm}\\
{\it Step 1. A priori $H^1$ estimate.}
It is easy to see that $y=y_1-y_2$ satisfies
\begin{equation}\label{y1-y2}\left\{
  \begin{array}{ll}
    -\nu\Delta  \mathbf y+ \mathbf{curl}\, \sigma(\mathbf y)\times  \mathbf y_2+
	\mathbf{curl}\, \sigma(\mathbf y_1)\times  \mathbf y+\nabla \pi= \mathbf u&\quad\mbox{in} \ \Omega,\vspace{2mm}\\
	 \mathrm{div} \,  \mathbf y=0&\quad\mbox{in} \ \Omega,\vspace{2mm}\\
   \mathbf y=0&\quad\mbox{on}\ \Gamma,
  \end{array}
\right.\end{equation}
where $u=u_1-u_2$. By setting $\phi=y$ in the corresponding variational formulation, we obtain
	$$\nu\left|y\right|_{H^1}^2+
	 \left( \mathbf{curl}\, \sigma(\mathbf y)\times \mathbf  y_2,  \mathbf y\right)=\left( u, y\right).$$
Due to Lemma \ref{rm2}, (\ref{state_est1}) and
 (\ref{state_est3}), it follows that
	$$\begin{array}{ll}\nu\left|y\right|_{H^1}^2
	&\leq \left\| u\|_2\| y\right\|_2+
\left(S_4^2\left|y_2\right|_{H^1}+\kappa
	\alpha\left\|y_2\right\|_{H^3}\right)
	\left|y\right|_{H^1}^2\vspace{2mm}\\
	&\leq S_2\left\|u\right\|_2\left|y\right|_{H^1} +
	\tfrac{\bar\kappa}{\nu} \left(\left\|u_2\right\|_2
	+\alpha \left\|{\rm curl}\, u_2\right\|_2\right)\left|y\right|_{H^1}^2
	\end{array}$$
 and (\ref{lipschitz_est_H1}) holds. \vspace{2mm}\\
{\it Step 2. A priori $H^3$ estimate.} To prove (\ref{lipschitz_est_V2}), let us first recall that if ${\rm curl}\,u_2\in H^1(\Omega)$, then ${\rm curl}\, \sigma(y_2)\in H^1(\Omega)$ (cf. Proposition \ref{uniqueness_state}). Using (\ref{transport_state}), we can see that $\tau={\rm curl}\, \sigma(y_2)-{\rm curl}\, \sigma(y_1)={\rm curl}\, \sigma(y)$ is the solution of the following transport equation
	$$\tau+\tfrac{\alpha}{\nu} \,y_1\cdot \nabla \tau
	+\tfrac{\alpha}{\nu} \, y\cdot 
	\nabla \left({\rm curl}\, \sigma(y_2)\right)
	={\rm curl}\, y
	+\tfrac{\alpha}{\nu}\,{\rm curl} \, u$$
and satisfies
	$$\left\|\tau\right\|_2^2+\tfrac{\alpha}{\nu} 
	\left( y\cdot 
	\nabla \left({\rm curl}\, \sigma(y_2)\right),\tau\right)=\left({\rm curl}\, y
	+\tfrac{\alpha}{\nu}\,{\rm curl} \, u,\tau\right).$$
By taking into account Lemma \ref{infty_V2} and using the Young inequality, we obtain
	$$\begin{array}{ll}
	\left\|\tau\right\|_2&\leq \left\|{\rm curl}\,y\right\|_2
	+\tfrac{\alpha}{\nu}
	\left\|{\rm curl} \, u\right\|_2+\tfrac{\alpha}{\nu}
	\left\|y\cdot 
	\nabla \left({\rm curl}\, \sigma(y_2)\right)
	\right\|_2\vspace{2mm}\\
	&\leq \left\|{\rm curl}\,y\right\|_2
	+\tfrac{\alpha}{\nu}
	\left\|{\rm curl} \, u\right\|_2+\tfrac{\alpha}{\nu} \left\|y\right\|_\infty \left|{\rm curl}\, \sigma(y_2)\right|_{H^1} \vspace{2mm}\\
	&\leq \left\|{\rm curl}\,y\right\|_2
	+\tfrac{\alpha}{\nu}
	\left\|{\rm curl} \, u\right\|_2+\tfrac{c\alpha^\frac{2}{3}}{\nu}\left|y\right|_{H^1}^{\frac{2}{3}}
	\left|{\rm curl}\, \sigma(y_2)\right|_{H^1} 
	\left\|\tau\right\|_2^{\frac{1}{3}}\vspace{2mm}\\
	&\leq \left\|{\rm curl}\,y\right\|_2
	+\tfrac{\alpha}{\nu}
	\left\|{\rm curl} \, u\right\|_2+
	\tfrac{1}{2} \left\|\tau\right\|_2+
	\tfrac{c\alpha}{\nu^\frac{3}{2}}
	\left|y\right|_{H^1}
	\left|{\rm curl}\, \sigma(y_2)\right|_{H^1}^\frac{3}{2}.\end{array}$$
Therefore,
	$$\left\|\tau\right\|_2
	\leq \tfrac{2\alpha}{\nu}\left\|{\rm curl}\, u\right\|_2
	+\kappa\left(\tfrac{\alpha}{\nu^\frac{3}{2}}
	\left|{\rm curl}\, \sigma(y_2)\right|_{H^1}^\frac{3}{2}+
	1\right)\left|y\right|_{H^1}$$
which gives the result.$\hfill\Box$

%%%%%%%%%%%%%%%%%%%%%%%%%%%%%%%%%%%%%%%%%%%%%%%%%%%%%%%
%%%%%%%%%%%%%%%%%%%%%%%%%%%%%%%%%%%%%%%%%%%%%%%%%%%%%%%
\subsection{G\^ateaux differentiability}
%\setcounter{equation}{0}
%%%%%%%%%%%%%%%%%%%%%%%%%%%%%%%%%%%%%%%%%%%%%%%%%%%%%%%
%%%%%%%%%%%%%%%%%%%%%%%%%%%%%%%%%%%%%%%%%%%%%%%%%%%%%%%
At this stage, we are able to study the G\^ateaux-differentiability of the control-to-state mapping. 
%%%%%%%%%%%%%%%%%%%%%%%%%%%%%%%%%%%%%%%%%%%%%%%%%%%%%%%
\begin{proposition}\label{taylor}
Let $u, w\in H(\mathrm{curl};\Omega)$ and assume in addition that 
${\rm curl}\, u\in H^1(\Omega)$ satisfies condition $(\ref{uniqueness_condition})$. For 
$0<\rho<1$, set  $ u_\rho= u+\rho  w$, and let $y$ and $ y_{\rho}$ be solutions of 
$(\ref{equation_etat})$ corresponding to $u$ and $u_\rho$, respectively. Then we have
	$$ y_{\rho}= y+\rho 
	 z+\rho r_\rho \qquad 
	\mbox{with} \ \lim_{\rho\rightarrow 0}
	\left|r_\rho\right|_{H^1}=0,$$
%and
%	$$J\left(u_\rho,y_\rho\right)=J\left(u,y\right)+\rho 
%	\left(\left(z,y-y_d\right)+\lambda(u,w)\right)
%	+o(\rho),$$
where $ z\in  V_2$ is a solution of $(\ref{linearized})$ corresponding to $(y,w)$.
\end{proposition}
%%%%%%%%%%%%%%%%%%%%%%%%%%%%%%%%%%%%%%%%%%%%%%%%%%%%%%%
{\bf Proof.} Easy calculation shows that $z_{\rho}=\frac{ y_{\rho}- y}{\rho}$ satisfies
	$$-\nu \Delta  \mathbf z_\rho+ \mathbf{curl}\, 
	\sigma\left(\mathbf z_\rho\right)\times \mathbf y+ \mathbf{curl}\,\sigma\left(\mathbf y_\rho\right)\times \mathbf z_\rho+\nabla\pi_\rho= \mathbf w.$$ 
Let $ z\in V_2$ be the solution of (\ref{linearized}).	
Then $r_\rho= z_\rho- z$ satisfies
	$$-\nu \Delta \mathbf r_\rho+ 
	\mathbf{curl}\,\sigma\left(\mathbf r_\rho\right)\times \mathbf y
	+ \mathbf{curl}\,\sigma\left(\mathbf y_\rho\right)\times
	 \mathbf r_\rho
	+ \mathbf{curl}\,\sigma\left(\mathbf y_\rho-\mathbf y\right)
	\times \mathbf z
	+\nabla\left(\pi_\rho-\pi\right)=0.$$
Multiplying this equation by $r_\rho$, we obtain	
	\begin{equation}\label{energy_r}\nu \left|r_\rho\right|_{H^1}^2+
	\left(\mathbf{curl}\,\sigma\left(\mathbf r_\rho\right)\times  \mathbf y+ \mathbf{curl}\,\sigma\left(\mathbf y_\rho\right)
	\times \mathbf r_\rho+\mathbf{curl}\,\sigma\left(\mathbf y_\rho-\mathbf y\right)
	\times \mathbf z,\mathbf r_\rho\right)=0.\end{equation}
It is easy to verify that
	\begin{equation}\label{energy_r2}\left( \mathbf{curl}\,\sigma\left(\mathbf y_\rho\right)
	\times \mathbf r_\rho,\mathbf r_\rho\right)=
	b\left(r_\rho,r_\rho,\sigma\left(y_\rho\right)\right)-
	b\left(r_\rho,r_\rho,\sigma\left(y_\rho\right)\right)=0.
	\end{equation}
Moreover, by taking into account Lemma \ref{rm2} and estimates (\ref{state_est1})-(\ref{state_est3}), we get
	\begin{equation}\label{energy_r3}\left|\left(\mathbf{curl}\,
	\sigma\left(\mathbf r_\rho\right)\times  \mathbf y,\mathbf 
	r_\rho\right)\right|
	\leq \tfrac{\bar\kappa}{\nu}\left(\left\| u\right\|_2+\alpha
	\left\| \mathrm{curl} \, u\right\|_2\right)
	\left|r_\rho\right|_{H^1}^2.
	\end{equation}
Combining (\ref{energy_r})-(\ref{energy_r3}), we deduce that
	$$\begin{array}{ll}
	\left(1-\tfrac{ \bar\kappa}{\nu^2}\left(\left\| u\right\|_2+\alpha
	\left\| \mathrm{curl} \, u\right\|_2\right)\right)
	\left|r_\rho\right|_{H^1}^2&\leq \tfrac{1}{\nu}
	\left|\left(\mathbf{curl}\,\sigma\left(\mathbf y_\rho-\mathbf y\right)
	\times \mathbf z,\mathbf r_\rho\right)\right|\vspace{2mm}\\
	&\leq \tfrac{1}{\nu}
	\left|y_\rho-y\right|_{V_2}
	 \left\|z\right\|_\infty
	\left\|r_\rho\right\|_2\vspace{2mm}\\
	&\leq \tfrac{S_2}{\nu}
	\left|y_\rho-y\right|_{V_2}
	 \left\|z\right\|_\infty
	\left|r_\rho\right|_{H^1}
	\end{array}$$
and thus
	$$\left(1-\tfrac{ \bar\kappa}{\nu^2}
	\left(\left\| u\right\|_2+\alpha
	\left\| \mathrm{curl} \, u\right\|_2\right)\right)
	\left|r_\rho\right|_{H^1}\leq \tfrac{S_2}{\nu}
	\left|y_\rho-y\right|_{V_2}
	 \left\|z\right\|_\infty.$$
The conclusion follows by observing that the term on the right-hand side of the previous inequality tends to zero when $\rho$ tends to zero. Indeed, due to (\ref{lipschitz_est_H1}) and (\ref{lipschitz_est_V2}),
 we have
	$$\begin{array}{ll}
	\left|y_\rho-y\right|_{V_2}&
	\leq \tfrac{2\alpha}{\nu}\left\|{\rm curl}\, (u_\rho-u)\right\|_2
	+\kappa\left(\tfrac{\alpha}{\nu^\frac{3}{2}}
	\left|{\rm curl}\, \sigma(y)\right|_{H^1}^\frac{3}{2}+
	1\right)\left|y_\rho-y\right|_{H^1}\vspace{1mm}\\
	&\leq \left(\tfrac{2\alpha}{\nu}\left\|{\rm curl}\, w\right\|_2
	+\kappa\left(\tfrac{\alpha}{\nu^\frac{3}{2}}
	\left|{\rm curl}\, \sigma(y)\right|_{H^1}^\frac{3}{2}+
	1\right)  \tfrac{S_2\nu}{\nu^2-\bar\kappa
	\left(\left\| u\right\|_2+\alpha
	\left\| \mathrm{curl} \, u\right\|_2
	\right)} \left\|w\right\|_2\right)\rho\vspace{2mm}\\
	&\longrightarrow 0 \qquad \mbox{when} \ \rho \rightarrow 0\end{array}$$
and the claimed result is proved. $\hfill\Box$
%%%%%%%%%%%%%%%%%%%%%%%%%%%%%%%%%%%%%%%%%%%%%%%%%%%%%%%
%%%%%%%%%%%%%%%%%%%%%%%%%%%%%%%%%%%%%%%%%%%%%%%%%%%%%%%
\section{Adjoint equation}
\setcounter{equation}{0}
%%%%%%%%%%%%%%%%%%%%%%%%%%%%%%%%%%%%%%%%%%%%%%%%%%%%%%%
%%%%%%%%%%%%%%%%%%%%%%%%%%%%%%%%%%%%%%%%%%%%%%%%%%%%%%%
Let  $u\in H( \mathrm{curl};\Omega)$ and let $ y\in V_2$ be a corresponding solution of (\ref{equation_etat}). The aim of this section is to study the existence of a weak solution for the adjoint state equation defined by
\begin{equation}\label{adjoint}\left\{
  \begin{array}{ll}
    -\nu\Delta  \mathbf p-  \mathbf{curl}\,\sigma(\mathbf y)\times  \mathbf p +\mathbf{curl}\left(\sigma\left(\mathbf y\times \mathbf p\right)\right)+\nabla \pi=\mathbf f&\quad\mbox{in} \ \Omega,\vspace{2mm}\\
	 \mathrm{div} \,  \mathbf p=0&\quad\mbox{in} \ \Omega,\vspace{2mm}\\
    \mathbf p=0 &\quad\mbox{on}\ \Gamma,
  \end{array}
\right .\end{equation}
where $f\in L^2(\Omega)$.
The two identities in Lemma \ref{non_lin_curl} motivates the following variational formulation.
%%%%%%%%%%%%%%%%%%%%%%%%%%%%%%%%%%%%%%%%%%%%%%%%%%%%%%%	
\begin{definition} \label{var_adj_2}A function $p\in  V$ is a weak solution of $(\ref{adjoint})$ if
	\begin{equation}\label{form_var_lin_adj}
	\nu\left(\nabla p,\nabla \phi\right)+
	b\left( p,\phi, \sigma(y) \right)-b\left( \phi,p, \sigma(y) \right)+
	b\left(p, y, \sigma(\phi) \right)-b\left( y,p, \sigma(\phi) \right)=
	\left(f, \phi\right)\end{equation}
for all $\phi\in V \cap H^2(\Omega)$.
\end{definition}
This formulation allows us to relate the adjoint state to the solution of the linearized equation and is particularly suited to derive the necessary optimality conditions. \vspace{1mm}\\
As will be seen below, existence of a Galerkin approximate solution
can be established by taking into account the formulation stated in Definition \ref{var_adj_2}. A corresponding a priori $H^1$ estimate can be derived and is sufficient to pass to the limit and prove the existence of a weak solution for the adjoint equation. Establishing a uniqueness result is much more challenging and requires higher regularity of the solutions. In this context, the observations raised in Section  \ref{lin_section}, concerning the most appropriate choice for the Galerkin basis, would similarly apply but deriving a $V_2$ estimate for the  approximate solution of (\ref{adjoint}) is far more difficult than in the case of the linearized equation. In order to illustrate our point, we can adapt the decomposition method and easily see that the term ${\rm curl}\, \sigma(p)$ should (formally) satisfy
	$${\rm curl}\, \sigma(p)-\tfrac{\alpha}{\nu}\,
	p\cdot \nabla\left({\rm curl}\, \sigma(y)\right)+
	\tfrac{\alpha}{\nu}\, {\rm curl}
	\left({\rm curl}\left(\sigma\left( y\times p\right)\right)\right)=
	\tfrac{\alpha}{\nu}\,{\rm curl}\, f+{\rm curl}\, p.
	$$
Solving this equation is not an easy task: in addition to high order derivatives of $p$ that we need to manage, the coefficients in ${\rm curl}\left({\rm curl}\left(\sigma\left( y\times \cdot\right)\right)\right)$ also involve high order derivatives of the state variable $y$. Following the ideas developed in Section 
\ref{existence_section}, we may prove that for every integer $k\geq 0$, if $\Gamma$ is of class $C^{k+2,1}$ then the semi-norm
 $\left|{\rm curl}\, \sigma(\cdot)\right|_{H^k}$ is equivalent to the norm
	$\left\|\cdot\right\|_{H^{k+2}}$. Recalling that  ${\rm curl}\, \sigma(y)$ satisfies (\ref{transport_state}) and in view of the classical regularity results for transport equations (generally based on fixed point arguments) the high order derivatives of the state variable are well defined if we assume that the control is accordingly regular and if we impose an additional restriction on the size of $y$ (and consequently on the corresponding control). Unlike the linearized equation, where the condition on the size of the data is set on the natural space $H({\rm curl};\Omega)$ and also guarantees uniqueness of the solution for the state equation and G\^ateaux differentiability of the control-to-state variable, the condition we need to impose here is set on higher-order Sobolev spaces and is much more restrictive.  \vspace{1mm}\\
%However, in the face of the difficulties encountered in establishing $H^2$ or $H^3$ a priori estimates for the approximate adjoint state, these considerations do not appear to be necessarily relevant.
An other aspect reinforces the idea that the effort in obtaining such   regularity results for the adjoint state may not be necessarily compensated. Keeping in mind that our objective is to derive first-order optimality conditions and that the natural framework for the controls is $H({\rm curl};\Omega)$, we should not require a priori additional regularity on this variable (and on the corresponding state). On the other hand, 
the results obtained in the previous sections concerning the solvability of the linearized state equation and the differentiability of the 
control-to-state mapping are only available in the case of regular data. To overcome this difficulty, our idea is to consider an approximate optimal control problem governed by a state equation involving regularized controls. The results stated in Sections
\ref{state_equation_section} and \ref{sec_lip} are then valid and we can derive the corresponding approximate optimality conditions. In order to pass to the limit, when the regularization parameter tends to zero, we only need a uniform estimate for the {\it regularized} adjoint state in $V$. 
%%%%%%%%%%%%%%%%%%%%%%%%%%%%%%%%%%%%%%%%%%%%%%%%%%%%%%%%
 \begin{proposition} \label{ex_uniq_adj}Let $ u\in H(\mathrm{curl};\Omega)$ satisfying condition $(\ref{uniqueness_condition})$ and let $ y\in  V_2$ be the corresponding solution of $(\ref{equation_etat})$. Then equation $(\ref{adjoint})$ admits at least a weak  solution $p\in  V$. Moreover, the following etimate holds
	\begin{equation}\label{est_H1_p}	
	\left(1-\tfrac{ \bar\kappa}{\nu^2}\left(\left\| u\right\|_2+\alpha
	\left\| \mathrm{curl} \, u\right\|_2\right)\right)\left|p\right|_{H^1}
	\leq \tfrac{S_2}{\nu} \left\|f\right\|_2.\end{equation}
\end{proposition}
%%%%%%%%%%%%%%%%%%%%%%%%%%%%%%%%%%%%%%%%%%%%%%%%%%%%%%%%
{\bf Proof.} We first establish the existence of an approximate solution and derive a corresponding apriori estimate in $H^1(\Omega)$. We next  pass to the limit and prove our statement.
\vspace{2mm}\\
{\it Step 1. Existence of an approximate solution and a priori $H^1$ estimate.} Consider the approximate problem defined by
	\begin{equation}\label{faedo_galerkin_adj}
	\left\{\begin{array}{ll}\mbox{Find} \ 
	 p_m=\displaystyle\sum_{j=1}^m \zeta_{j}
	  e_j \ \mbox{solution}, \ \mbox{for} \ 1\leq j\leq m, \ 
	\mbox{of}\vspace{3mm}\\
	\nu\left(\nabla p_m,\nabla e_j\right)-\left( \mathbf{curl}\,\sigma(\mathbf y)\times  \mathbf p_m, e_j\right)+\left(\mathbf{curl}\left(\sigma\left(\mathbf y\times \mathbf p_m\right)\right), e_j\right)=
	\left(f,  e_j\right).\end{array}\right.
	\end{equation}
where $( e_j)_j\subset H^4(\Omega)$ is the set of the eigenfunctions, solutions of (\ref{eigenfunc}). 
Due to Lemma \ref{non_lin_curl},  we have
	\begin{align}\label{form_adj_vw}
	-\left( \mathbf{curl}\,\sigma(\mathbf y)\times  \mathbf p_m, e_j\right)
	+\left(\mathbf{curl}\left(\sigma\left(\mathbf y\times \mathbf p_m\right)\right), e_j\right)&
	=b\left( p_m, e_j, \sigma(y) \right)-b\left( e_j,p_m, \sigma(y) \right)\nonumber\\
	&+
	b\left(p_m, y, \sigma( e_j) \right)-b\left( y,p_m, \sigma( e_j) \right).\end{align}
Let then $m$ be fixed and consider $Q: \ \mathbb{R}^m \longrightarrow \mathbb{R}^m$ defined by
	$$\begin{array}{ll}\left(Q \zeta\right)_i=&\nu\left(\nabla p_m,\nabla  e_j\right)+
	b\left( p_m, e_j, \sigma(y) \right)-b\left( e_j,p_m, \sigma(y) \right)\vspace{2mm}\\
	&+
	b\left(p_m, y, \sigma( e_j) \right)-b\left( y,p_m, \sigma( e_j) \right)-
	\left(f,   e_j\right),	\end{array}$$
where $ p_m=\sum_{i=1}^m\zeta_i   e_i$. The mapping $Q$ is obviously continuous. Let us prove that $Q(\zeta)\cdot \zeta>0$ if $|\zeta|$ is sufficiently large. Arguing as in the proof of  Proposition \ref{ex_uniq_lin}, we may prove that 
	\begin{align}\label{Q_zeta}{Q}(\zeta)\cdot \zeta&=
	\nu \left|p_m\right|_{H^1}^2+b\left(p_m, y, \sigma(p_m) \right)-b\left( y,p_m, \sigma(p_m) \right)-
	\left(f, p_m\right)\\
	&\geq \left(\nu-\tfrac{\bar\kappa}{\nu} \left(\left\| u\right\|_2+\alpha
	\left\| \mathrm{curl} \, u\right\|_2\right)
	\right)\left|\zeta\right|^2-
	\left\|f\right\|_2\left|\zeta\right|\nonumber\\
	&\longrightarrow +\infty \qquad \mbox{when} \ 
	|\zeta|\rightarrow +\infty.\nonumber\end{align}
Due to the Brouwer theorem, we deduce that there exists $ \zeta^\ast\in \mathbb{R}^k$ such that ${Q}\left( \zeta^\ast\right)=0$ and thus $ p_m=\sum_{i=1}^k\zeta_i^\ast   e_i$ is a solution of problem (\ref{faedo_galerkin_adj}). Taking into account $(\ref{Q_zeta})$ and  Lemma \ref{rm2}, we deduce that
	$$\begin{array}{ll}\nu \left|p_m\right|_{H^1}^2&=
	\left( \mathbf{curl}\,\sigma(\mathbf p_m)\times  \mathbf y, \mathbf  p_m\right)-
	\left(f, p_m\right)\vspace{2mm}\\
	&\leq\tfrac{ \bar\kappa}{\nu}\left(\left\| u\right\|_2+\alpha
	\left\| \mathrm{curl} \, u\right\|_2\right)
	\left|p_m\right|_{H^1}^2+
	S_2\left\|f\right\|_2
	 \left|p_m\right|_{H^1}
	\end{array}$$
which gives
	\begin{equation}\label{est_H1_pm}\left(
	1-\tfrac{ \bar\kappa}{\nu^2}\left(\left\| u\right\|_2+\alpha
	\left\| \mathrm{curl} \, u\right\|_2\right)\right)\left\|\nabla p_m\right\|_2
	\leq \tfrac{S_2}{\nu} \left\|f\right\|_2.\end{equation}
{\it Step 2. Passing to the limit.} It remains to pass to the limit with respect to $m$. From estimate (\ref{est_H1_pm}), it follows that if 
$\bar u$ satisfies condition $(\ref{uniqueness_condition})$, then there exists a subsequence, still indexed by $m$, and function $ p\in  V_2$ such that
	$$ p_m \longrightarrow  p \qquad \mbox{weakly in} \  V.$$
By taking into account (\ref{form_adj_vw}) and passing to the limit in (\ref{faedo_galerkin_adj}), we obtain for every $j\geq 1$
	$$\nu\left(\nabla p,\nabla e_j\right)+b\left( p, e_j, \sigma(y) \right)-b\left( e_j,p, \sigma(y) \right)+
	b\left(p, y, \sigma( e_j) \right)-b\left( y,p, \sigma( e_j) \right))=
	\left(f,   e_j\right)	$$
and by density we prove that $ p$ satisfies the variational formulation (\ref{form_var_lin_adj}). Moreover, $p$ satisfies (\ref{est_H1_p}). $\hfill\Box$

%%%%%%%%%%%%%%%%%%%%%%%%%%%%%%%%%%%%%%%%%%%%%%%%%%%%%%%
%%%%%%%%%%%%%%%%%%%%%%%%%%%%%%%%%%%%%%%%%%%%%%%%%%%%%%%
\section{Proof of the main results}
\setcounter{equation}{0}
Unless necessary, and in order to simplify the redaction, the index $\alpha$ will be dropped.
%%%%%%%%%%%%%%%%%%%%%%%%%%%%%%%%%%%%%%%%%%%%%%%%%%%%%%%
%%%%%%%%%%%%%%%%%%%%%%%%%%%%%%%%%%%%%%%%%%%%%%%%%%%%%%%
\subsection{Proof of the existence of an optimal control for $(P_\alpha)$}
%%%%%%%%%%%%%%%%%%%%%%%%%%%%%%%%%%%%%%%%%%%%%%%%%%%%%%%
We first prove Theorem \ref{main_existence}. Let $(u_k,y_k)_k\subset U_{ad}\times V_2$ be a minimizing sequence. Since $(u_k)_k$ is uniformly bounded in the closed convex set $U_{ad}$, we may extract a subsequence, still indexed by $k$, weakly convergent to some $u\in U_{ad}$. Applying Proposition \ref{continuity_H3} with $U=U_{ad}$, it follows that $(y_k)_k$ converges to $y$, solution of (\ref{equation_etat}) corresponding to $u$, in $H^2(\Omega)$. The convexity and continuity of $J$ imply the lower semicontinuity of $J$ in the weak topology and 
	$$J(u,y)\leq \liminf_k J(u_k,y_k)=\inf(P_\alpha),$$
showing that $(u,y)$ is a solution for $(P_\alpha)$.$\hfill\Box$
%%%%%%%%%%%%%%%%%%%%%%%%%%%%%%%%%%%%%%%%%%%%%%%%%%%%%%
\subsection{\bf Proof of the necessary optimality conditions 
 for $(P_\alpha)$}
%%%%%%%%%%%%%%%%%%%%%%%%%%%%%%%%%%%%%%%%%%%%%%%%%%%%%%
\subsubsection{Approximate optimal control problem}
%%%%%%%%%%%%%%%%%%%%%%%%%%%%%%%%%%%%%%%%%%%%%%%%%%%%%%
For $\varepsilon>0$, we denote by $\varrho_\varepsilon$ a Friedrichs mollifier, i.e. the convolution operator defined by 
	$$\varrho_\varepsilon(u)=\varrho_\varepsilon \ast u,$$
 where $\varrho_\varepsilon(x)=\varepsilon^{-2} \varrho\left(\tfrac{x}{\varepsilon}\right)$ and $\varrho$ is a positive radial compactly supported smooth function whose integral is equal $1$. Let us recall some  usefull properties of these mollifiers.
	\begin{enumerate}
	\item $\varrho_\varepsilon$ is selfadjoint for the $L^2$ scalar product, i.e.
	$$
	\left(\varrho_\varepsilon(u),v\right)=\left(u,\varrho_\varepsilon(v)\right)
	 \qquad \mbox{for all} \ u,v\in L^2(\Omega).$$
	\item $\varrho_\varepsilon$ commutes with derivatives.
	\item For $m\in \mathbb{N}$ and $u\in H^m(\Omega)$, we have
	$$
	\left\|\varrho_\varepsilon(u)\right\|_{H^m}\leq 
	\left\|u\right\|_{H^m} \quad \mbox{and} \quad 
	\lim_{\varepsilon \rightarrow 0^+}
	\left\|u-\varrho_\varepsilon(u)\right\|_{H^m}=0.$$
	\end{enumerate}
Due to the first and second properties, we have 
	\begin{align}\label{auto_adj_molli}
	\left(\varrho_\varepsilon \left(u\right),v
	\right)_{H({\rm curl};\Omega)}&=
	\left(\varrho_\varepsilon \left(u\right),v\right)+
	\left({\rm curl}\, \varrho_\varepsilon \left(u\right),
	{\rm curl}\,v\right)\nonumber\\
	&=\left(\varrho_\varepsilon \left(u\right),v\right)+
	\left(\varrho_\varepsilon \left({\rm curl}\, u\right),
	{\rm curl}\,v\right)\nonumber\\
	&=\left(u,\varrho_\varepsilon \left(v\right)\right)+
	\left({\rm curl}\, u,
	\varrho_\varepsilon \left({\rm curl}\,v\right)\right)\nonumber\\
	&=\left(u,\varrho_\varepsilon \left(v\right)\right)+
	\left({\rm curl}\, u,
	{\rm curl}\, \varrho_\varepsilon \left(v\right)\right)\nonumber\\
	&=\left(u,\varrho_\varepsilon \left(v\right)
	\right)_{H({\rm curl};\Omega)}
	\qquad \mbox{for all} \ u,v\in H({\rm curl};\Omega).\end{align}
Let $(\bar u,\bar y)$ be a solution of $(P_\alpha)$ with $\bar u$ satisfying condition (\ref{uniqueness_condition}) and consider the control problem $(P_\alpha^\varepsilon)$ defined in Section \ref{introduction}. We first prove the existence of an optimal control for $(P_\alpha^\varepsilon)$. The proof combines the standard arguments already used to establish Theorem \ref{main_existence} with the properties of the mollifiers.
%%%%%%%%%%%%%%%%%%%%%%%%%%%%%%%%%%%%%%%%%%%%%%%%%%%%%%
\begin{proposition} Assume that $U_{ad}$ is bounded in $H(\mathrm{curl};\Omega)$. Then problem $(P_\alpha^\varepsilon)$ admits a solution.
\end{proposition}
%%%%%%%%%%%%%%%%%%%%%%%%%%%%%%%%%%%%%%%%%%%%%%%%%%%%%%
{\bf Proof.} Let $(u_{k}^\varepsilon,y_{k}^\varepsilon)_k\subset U_{ad}\times V_2$ be a minimizing sequence for $(P_\alpha^\varepsilon)$. Then there exists a subsequence, still indexed by $k$, and $u^\varepsilon\in U_{ad}$ such that $(u_k^\varepsilon)_k$ weakly converges to $u^\varepsilon$ in $H(\mathrm{curl};\Omega)$. Since $\varrho_\varepsilon$ is linear and continuous, it is weakly continuous and thus $\left(\varrho_\varepsilon \left(u_k^\varepsilon\right)\right)_k$ weakly converges to $\varrho_\varepsilon \left(u^\varepsilon\right)$ 
in $H({\rm curl};\Omega)$.
By taking into account Proposition \ref{continuity_H3} (with $U=\varrho_\varepsilon(U_{ad})$), we deduce that $\left(y_{k}^\varepsilon\right)_{k}$ converges in $H^2(\Omega)$ to
$y^\varepsilon$, a solution of (\ref{equation_etat}) corresponding to $\varrho_\varepsilon \left(u^\varepsilon\right)$. This implies that $(y^\varepsilon,u^\varepsilon)$ is admissible for $(P_\alpha^\varepsilon)$ and by using the convexity and continuity of $I$, we obtain 
	$$I(u^\varepsilon,y^\varepsilon)\leq \liminf_k I(u_k^\varepsilon,y_k^\varepsilon)
	=\inf(P_\alpha^\varepsilon),$$
showing that $(u^\varepsilon,y^\varepsilon)$ is a solution for $(P_\alpha^\varepsilon)$.$\hfill\Box$\vspace{1mm}\\
%%%%%%%%%%%%%%%%%%%%%%%%%%%%%%%%%%%%%%%%%%%%%%%%%%%%%%
The next result deals with the necessary optimality conditions for the approximate problem $(P_\alpha^{\varepsilon})$.
%%%%%%%%%%%%%%%%%%%%%%%%%%%%%%%%%%%%%%%%%
\begin{proposition} \label{aprox_optimal_cond} Let $(\bar u^\varepsilon,\bar y^\varepsilon)$ be a solution of $(P_\alpha^\varepsilon)$ and assume that $\bar u^\varepsilon$ satisfies condition $(\ref{uniqueness_condition})$. Then there exists $\bar{p}^\varepsilon\in V$ such that
	\begin{equation}\label{state_regularized}
	\left\{ \begin{array}{ll}-\nu \Delta  \mathbf{\bar y}^\varepsilon+
	\mathbf{curl}\,\sigma\left(\mathbf{\bar y}^\varepsilon	\right)\times  \mathbf{\bar y}^\varepsilon+\nabla \pi^\varepsilon= \varrho_\varepsilon(\mathbf{\bar u}^\varepsilon)& \mbox{in} \ \Omega,\vspace{2mm} \\
             \mathrm{div} \,  \mathbf{\bar y}^\varepsilon=0& \mbox{in} \ \Omega,\vspace{2mm}\\
	\mathbf{\bar y}^\varepsilon=0&
	\mbox{on} \ \Gamma,\end{array}\right.
	\end{equation}
	\begin{equation}\label{adjoint_regularized}\left\{ \begin{array}{ll}-\nu \Delta  \mathbf{\bar p}^\varepsilon-
	\mathbf{curl}\,\sigma\left(\mathbf{\bar y}^\varepsilon\right)
	\times  \mathbf{\bar p}^\varepsilon+\mathbf{curl}\,\sigma\left(\mathbf{\bar y}^\varepsilon\times \mathbf{\bar p}^\varepsilon\right)
	+\nabla \tilde\pi^\varepsilon= \mathbf{\bar y}^\varepsilon-\mathbf{y_d}& \mbox{in} \ \Omega,\vspace{2mm} \\
             \mathrm{div} \,  \mathbf{\bar p}^\varepsilon=0& \mbox{in} \ \Omega,\vspace{2mm}\\
	\mathbf{\bar p}^\varepsilon=0&
	\mbox{on} \ \Gamma,\end{array}\right.\end{equation}
	\begin{equation}\label{op_control}
	\left(\varrho_\varepsilon\left(\bar{p}^{\varepsilon}\right)+
	\lambda\bar{u}^{\varepsilon},
	v-\bar{u}^{\varepsilon}\right)+
	\left(\bar u^\varepsilon-\bar u,v-\bar u^\varepsilon\right)_{H({\rm curl};\Omega)}\geq 0\qquad \mbox{for all} \
	  v\in U_{ad}.\end{equation}
\end{proposition}
%%%%%%%%%%%%%%%%%%%%%%%%%%%%%%%%%%%%%%%%%%%%%%%%%%%%%%
{\bf Proof.} Let first notice that if $\bar{u}^{\varepsilon}$ satisfies condition $(\ref{uniqueness_condition})$, then $\varrho_\varepsilon\left(\bar{u}^{\varepsilon}\right)$ satisfies the same condition.
Taking into account Proposition \ref{uniqueness_state}, we deduce that (\ref{state_regularized}) admits a unique solution ${\bar y}^\varepsilon$ and that this solution belongs to $H^4(\Omega)$. Due to Proposition \ref{ex_uniq_lin}, it follows that for every $v\in H(\mathrm{curl};\Omega)$, the linearized equation
	\begin{equation}\label{linearized_regularized}
	\left\{ \begin{array}{ll}-\nu \Delta  \mathbf z+
	\mathbf{curl}\,\sigma\left(\mathbf z\right)\times \mathbf{\bar y}^\varepsilon+\mathbf{curl}\,\sigma\left(\mathbf{\bar y}^\varepsilon\right)\times  \mathbf z+\nabla \pi^\varepsilon= \varrho_\varepsilon(\mathbf v)& \mbox{in} \ \Omega,\vspace{2mm} \\
             \mathrm{div} \,  \mathbf z=0& \mbox{in} \ \Omega,\vspace{2mm}\\
	\mathbf z=0&
	\mbox{on} \ \Gamma,\end{array}\right.
	\end{equation}
admits a unique solution $\bar z^\varepsilon(v)\in H^3(\Omega)$. Moreover, due Proposition \ref{taylor} and  Proposition \ref{ex_uniq_adj},  the control-to-state mapping $u\mapsto y^\varepsilon(u)$ is G\^ateaux differentiable at $\bar u^\varepsilon$ and equation (\ref{adjoint_regularized}) admits at least a solution ${\bar p}^\varepsilon\in V$.\vspace{2mm}\\
 For $\rho\in ]0,1[$ and $v\in U_{ad}$, let $u_{\rho}^\varepsilon=\bar u^\varepsilon+\rho (v-\bar u^\varepsilon)$, $y_\rho^\varepsilon$ be the solution of (\ref{equation_etat}) corresponding to $\varrho_\varepsilon\left(u_{\rho}^\varepsilon\right)$ and $z_{\rho}^\varepsilon=\tfrac{y_{\rho}^\varepsilon-\bar y^\varepsilon}{\rho}$. Since $(\bar u^\varepsilon,\bar y^\varepsilon)$ is an optimal solution for $(P_\alpha^\varepsilon)$ and $(u_{\rho}^\varepsilon,y_{\rho}^\varepsilon)$ is admissible for this problem, we have 
	$$\displaystyle\lim_{\rho\rightarrow 0}
	\tfrac{I(u_{\rho}^\varepsilon,y_{\rho}^\varepsilon)-
	I(\bar u^\varepsilon,\bar y^\varepsilon)}{\rho}\geq 0$$
which yields
	\begin{equation}\label{zv}\left(\bar z^\varepsilon(v-\bar u^\varepsilon),\bar y^\varepsilon-y_d\right)+
	\lambda\left(\bar u^\varepsilon,v-\bar u^\varepsilon\right)+
	\left(\bar u^\varepsilon-\bar u,v-\bar u^\varepsilon\right)_{H({\rm curl};\Omega)}\geq 0.\end{equation}
 Setting $\phi=\bar z^\varepsilon(v-\bar u^\varepsilon)$ in the variational formulation 
(\ref{form_var_lin_adj}) corresponding to $\bar p^\varepsilon$ and taking into account the variational formulation (\ref{var_lin}), we obtain
\begin{align}\left(\bar{y}^\varepsilon-y_d,\bar z^\varepsilon(v-\bar u^\varepsilon)\right)
&=\nu\left(\nabla \bar p^\varepsilon,\nabla \bar z^\varepsilon(v-\bar u^\varepsilon)\right)+
	b\left(\bar p^\varepsilon,\bar z^\varepsilon(v-\bar u^\varepsilon), \sigma\left(\bar y^\varepsilon\right) \right)
	-b\left(\bar z^\varepsilon(v-\bar u^\varepsilon),{\bar p}^\varepsilon, \sigma\left(\bar y^\varepsilon\right) \right)\nonumber\\
& \ +
	b\left(\bar p^\varepsilon, \bar y^\varepsilon, \sigma\left(\bar z^\varepsilon(v-\bar u^\varepsilon)\right) \right)
	-b\left(\bar y^\varepsilon,\bar p^\varepsilon, \sigma\left(\bar z^\varepsilon(v-\bar u^\varepsilon)\right) \right)\nonumber\\
&\label{GF_op}=\left(\varrho_\varepsilon\left(v-\bar{u}^\varepsilon\right),\bar{p}^\varepsilon\right)=\left(v-\bar{u}^\varepsilon,\varrho_\varepsilon\left(\bar{p}^\varepsilon\right)\right).
	\end{align}
The result follows by combining (\ref{zv}) and (\ref{GF_op}).$\hfill\Box$
%%%%%%%%%%%%%%%%%%%%%%%%%%%%%%%%%%%%%%%%%%%%%%%%%%%%%%
\begin{proposition} \label{conv_con_ap}Let  $({\bar u}^\varepsilon,{\bar y}^\varepsilon)$ be a solution for $(P_\alpha^\varepsilon)$. There exists a subsequence $(\varepsilon_k)_k$ converging to zero, such that
	$$\lim_{k\rightarrow +\infty}\left\|{\bar u}^{\varepsilon_k}-\bar u\right\|_{H({\rm curl};\Omega)}=0$$
	$$\lim_{k\rightarrow +\infty}\left\|{\bar y}^{\varepsilon_k}-\bar y\right\|_{H^3}=0$$
	$$\lim_{k\rightarrow +\infty}I\left({\bar u}^{\varepsilon_k},{\bar y}^{\varepsilon_k}\right)=J(\bar u,\bar y)$$
\end{proposition}
%%%%%%%%%%%%%%%%%%%%%%%%%%%%%%%%%%%%%%%%%%%%%%%%%%%%%%
{\bf Proof.} Since $(\bar u^\varepsilon)_\varepsilon$ is bounded in $U_{ad}$, there exists a subsequence $(\varepsilon_k)_k$ converging to zero and $u\in U_{ad}$ such that $(\bar u^{\varepsilon_k})_k$ converges to $u$ weakly in $H(\mathrm{curl};\Omega)$. Due to (\ref{auto_adj_molli})
we have
	$$\begin{array}{ll}
	\displaystyle\lim_{k\rightarrow +\infty}\left(\varrho_{\varepsilon_k}
	\left(\bar u^{\varepsilon_k}\right),
	\phi\right)_{H({\rm curl};\Omega)}&=\displaystyle
	\lim_{k\rightarrow +\infty}\left(\bar u^{\varepsilon_k},
	\varrho_{\varepsilon_k}\left(\phi\right)\right)_{H({\rm curl};\Omega)}
	\vspace{2mm}\\
	&=\displaystyle\left(u,\phi\right)_{H({\rm curl};\Omega)}
	\qquad \mbox{for all} \ \phi\in H(\mathrm{curl};\Omega)\end{array}$$
implying that $\left(\varrho_{\varepsilon_k}
	\left(\bar u^{\varepsilon_k}\right)\right)_k$ also weakly converges to $u$ in $H(\mathrm{curl};\Omega)$. Due to Proposition 
\ref{continuity_H3}, we deduce that $(\bar y^{\varepsilon_k})_k$
 converges in $H^2(\Omega)$ to $y$, solution of (\ref{equation_etat}) corresponding to $u$. On the other hand, let $y^{\varepsilon_k}_{\bar u}$ be a solution of (\ref{equation_etat}) corresponding to $\varrho_{\epsilon_k}(\bar u)$. Since $(\varrho_{\epsilon_k}(\bar u))_k$ strongly converges to $\bar u$ in $H(\mathrm{curl};\Omega)$, it follows that $\left(y^{\varepsilon_k}_{\bar u}\right)_k$ converges to $\bar y$ in $H^3(\Omega)$.  Using the lower semicontinuity of $I$ and the admissibility of $(\bar u,y_{\bar u}^{\varepsilon_k})$ for $(P_\alpha^{\varepsilon_k})$, we obtain
	$$\begin{array}{ll}\displaystyle\tfrac{1}{2}\left\|y-y_d\right\|^2_2+\tfrac{\lambda}{2}
	\left\|u\right\|^2_2+\tfrac{1}{2}\left\|u-\bar{u}\right\|^2_{H({\rm curl};\Omega)}&\leq\displaystyle\liminf_{k} 
	I(\bar{u}^{\varepsilon_k},\bar y^{\varepsilon_k})\vspace{2mm}\\
	&\displaystyle\leq\limsup_{k} I(\bar{u}^{\varepsilon_k},\bar y^{\varepsilon_k})\vspace{1mm}\\
	&\leq\displaystyle\lim_{k} I(\bar u,y_{\bar u}^{\varepsilon_k})
	=\tfrac{1}{2}\left\|\bar y-y_d\right\|^2_2+\tfrac{\lambda}{2}\left\|\bar{u}\right\|^2_2\end{array}$$
and consequently $$J(u,y)+\tfrac{1}{2}\left\|u-\bar{u}\right\|^2_{H({\rm curl};\Omega)}\leq J(\bar{u},\bar y).$$
Since $(\bar{u},\bar y)$ is solution of $(P_\alpha)$, we have $J(\bar{u},\bar y)
\leq J(u,y)$ and thus $u=\bar{u}$. Recalling that $\bar u$ satisfies condition (\ref{uniqueness_condition}), we deduce that $y=\bar y$ and thus 
         $$ \lim_{k\rightarrow +\infty}I(\bar{u}^{\varepsilon_k},\bar y^{\varepsilon_k})=J(\bar{u},\bar y).$$
Finally, observing that
	$$\begin{array}{ll}\displaystyle\tfrac{1}{2}\limsup_{k}\left\|\bar{u}^{\varepsilon_k}-\bar{u}\right\|^2_{H({\rm curl};\Omega)}
	&=\displaystyle\limsup_{k}\left(I(\bar{u}^{\varepsilon_k},\bar y^{\varepsilon_k})-\tfrac{1}{2}
	\|\bar y^{\varepsilon_k}-y_d\|^2_2-\tfrac{\lambda}{2}\left\|\bar u^{\varepsilon_k}
	\right\|^2_2\right)\vspace{1mm}\\
	&\displaystyle\leq J(\bar{u},\bar y)-\tfrac{1}{2}\|\bar{y}-y_d\|^2_2-\tfrac{\lambda}{2}\liminf_{k}
	\|\bar u^{\varepsilon_k}\|^2_2\vspace{2mm}\\
	&\displaystyle=\tfrac{\lambda}{2}\left\|\bar{u}\right\|_2^2-\tfrac{\lambda}{2}
	\liminf_{k}\left\|\bar{u}^{\varepsilon_k}\right\|^2_2\leq 0\end{array}$$
we conclude that $(\bar{u}^{\varepsilon_k})_k$ converges to $\bar{u}$ strongly in $H({\rm curl};\Omega)$.$\hfill\Box$
%%%%%%%%%%%%%%%%%%%%%%%%%%%%%%%%%%%%%%%%%%%%%%%%%%%%%%
%%%%%%%%%%%%%%%%%%%%%%%%%%%%%%%%%%%%%%%%%%%%%%%%%%%%%%
\subsubsection{Proof of Theorem \ref{main_1}} Let $(\bar u^{\varepsilon_k},\bar y^{\varepsilon_k})$ be the solution of $(P^{\varepsilon_k}_\alpha)$ given in Proposition \ref{conv_con_ap}. Since $\bar u$ satisfies condition $(\ref{uniqueness_condition})$, we deduce that there exists $k_1\in \NN$ such that $\bar u^{\varepsilon_k}$ also satisfies condition $(\ref{uniqueness_condition})$ for every $k> k_1$. Due Proposition \ref{aprox_optimal_cond}, there exists  $\bar{p}^{\varepsilon_k}\in V$ such that
\begin{equation}\label{adjoint_regularized_k}\left\{ \begin{array}{ll}-\nu \Delta  \mathbf{\bar p}^{\varepsilon_k}-
	\mathbf{curl}\,\sigma\left(\mathbf{\bar y}^{\varepsilon_k}\right)
	\times  \mathbf{\bar p}^{\varepsilon_k}+\mathbf{curl}\,\sigma\left(\mathbf{\bar y}^{\varepsilon_k}\times \mathbf{\bar p}^{\varepsilon_k}\right)
	+\nabla \tilde\pi^{\varepsilon_k}= \mathbf{\bar y}^{\varepsilon_k}-\mathbf{y_d}& \mbox{in} \ \Omega,\vspace{2mm} \\
             \mathrm{div} \,  \mathbf{\bar p}^\varepsilon_k=0& \mbox{in} \ \Omega,\vspace{2mm}\\
	\mathbf{\bar p}^{\varepsilon_k}=0&
	\mbox{on} \ \Gamma,\end{array}\right.\end{equation}
	\begin{equation}\label{op_control_k}\left(\varrho_{\varepsilon_k}\left(\bar{p}^{\varepsilon_k}\right)+
	(\lambda+1)\bar{u}^{\varepsilon_k}-\bar{u},
	v-\bar{u}^{\varepsilon_k}\right)\geq 0\qquad \mbox{for all} \
	  v\in U_{ad}.\end{equation}
Moreover, due to (\ref{est_H1_p}),  we have the following estimate
	\begin{equation}\label{est_H1_p_k}	
	\left|p^{\varepsilon_k}\right|_{H^1}
	\leq \tfrac{S_2\nu}{\nu^2-\bar\kappa
	\left(\left\| u^{\varepsilon_k}\right\|_2+\alpha
	\left\| \mathrm{curl} \, u^{\varepsilon_k}\right\|_2\right)} \left\|\bar y^{\varepsilon_k}-y_d\right\|_2.\end{equation}
Using once again the strong convergence of $(u^{\varepsilon_k},y^{\varepsilon_k})_k$ in $H({\rm curl; \Omega})\times L^2(\Omega)$, we deduce that the sequence $(p^{\varepsilon_k})_k$
is bounded in $V$. There then exist a subsequence, still indexed by $k$, 
and $\bar p$ such that $(\bar p^{\varepsilon_k})_{k}$ weakly converges to $\bar p$ in $V$ and, by using compactness results on Sobolev spaces, $(\bar p^{\varepsilon_k})_{k}$ strongly converges to $\bar p$ in $L^{2}(\Omega)$. Therefore, by passing into the limit in the variational formulation corresponding to $\bar p^{\varepsilon_k}$:
	$$ \begin{array}{ll}
	&\nu\left(\nabla \bar p^{\varepsilon_k},\nabla \phi\right)+
	b\left(\bar p^{\varepsilon_k},\phi, 
	\sigma\left(\bar y^{\varepsilon_k}\right) \right)
	-b\left( \phi,\bar p^{\varepsilon_k}, \sigma\left(\bar y^{\varepsilon_k}\right) \right)+
	b\left(\bar p^{\varepsilon_k}, \bar y^{\varepsilon_k}, \sigma(\phi) \right)-b\left(\bar y^{\varepsilon_k},\bar p^{\varepsilon_k}, \sigma(\phi) \right)\vspace{2mm}\\
	&=
	\left(\bar y^{\varepsilon_k}-y_d, \phi\right)\end{array}$$
we obtain
	$$\nu\left(\nabla \bar p,\nabla \phi\right)+
	b\left(\bar p,\phi, 
	\sigma\left(\bar y\right) \right)
	-b\left( \phi,\bar p, \sigma\left(\bar y\right) \right)+
	b\left(\bar p, \bar y, \sigma(\phi) \right)-b\left(\bar y,\bar p, \sigma(\phi) \right)=
	\left(\bar y-y_d, \phi\right)$$
for all $\phi\in V \cap H^2(\Omega)$, that is $\bar p$ is a weak solution of (\ref{adj_opt_eq_alpha}). Finally, observing that
	$$\begin{array}{ll}\left\|\varrho_{\varepsilon_k}
	\left(\bar p^{\varepsilon_k}\right)-\bar p\right\|_2&\leq 
	\left\|\varrho_{\varepsilon_k}
	\left(\bar p^{\varepsilon_k}-\bar p\right)\right\|_2+
	\left\|\varrho_{\varepsilon_k}
	\left(\bar p\right)-\bar p\right\|_2\leq 
	\left\|\bar p^{\varepsilon_k}-\bar p\right\|_2+
	\left\|\varrho_{\varepsilon_k}
	\left(\bar p\right)-\bar p\right\|_2\vspace{2mm}\\
	&\longrightarrow 0 \quad \mbox{when} \ k\rightarrow \infty,\end{array}$$
we obtain (\ref{opt_control_alpha}) by passing into the limit in (\ref{op_control_k}).$\hfill\Box$

%%%%%%%%%%%%%%%%%%%%%%%%%%%%%%%%%%%%%%%%%%%%%%
\subsection{Asymptotic analysis when $\alpha$ tends to zero}
%%%%%%%%%%%%%%%%%%%%%%%%%%%%%%%%%%%%%%%%%%%%%%
 The proof of Theorem  \ref{assympt_1} is split into three steps. First, we prove that if $(u_\alpha,y_\alpha)$ is admissible for $(P_\alpha)$, then it  converges in the weak-$H(curl;\Omega)\times V$ topology to an admissible point $(u_0,y_0)$ for problem $(P_0)$. Next, we prove that if $(\bar u_\alpha,\bar y_\alpha)$ is an optimal solution of $(P_\alpha)$, then the limit point $(\bar u_0,\bar y_0)$ is an optimal solution of $(P_0)$ and the convergence of $\bar u_\alpha$ to $\bar u_0$ is strong in the topology of $L^2(\Omega)$. Finally, we pass to the limit in the adjoint equation and prove that the limit point $\bar p_0$ satisfies an adjoint equation and optimality condition associated with $(P_0)$.\vspace{2mm}\\
{\it Step 1. Convergence of admissible points.} Let $(u_\alpha,y_\alpha)$ be an admissible point for $(P_\alpha)$.
By taking into account (\ref{state_est1}) and (\ref{state_est2}), we have
	$$\left|y_\alpha\right|_2\leq \tfrac{S_2}{\nu} \left\|u_\alpha\right\|_2, $$
	$$ \left|y_\alpha\right|_{V_2}\leq \tfrac{1}{\nu}\left(
	S_2\left\|u_\alpha\right\|_2+\alpha \left\|\mathrm{curl}\,u_\alpha\right\|_2\right),$$
and thus $\left(y_\alpha\right)_\alpha$ and
$\left(\mathrm{curl}\, \sigma(y_\alpha)\right)_\alpha$ are bounded independently of $\alpha$. There then exists a subsequence, still indexed by $\alpha$, $u_0\in U_{ad}$, $y_0\in V$ and $\omega_0\in L^2(\Omega)$ such that
	$$u_\alpha\longrightarrow 
	u_0 \qquad
	 \mbox{weakly in} \ L^2(\Omega). $$
	$$y_\alpha \longrightarrow  y_0 \qquad
	 \mbox{weakly in} \ H^1(\Omega) \
	\mbox{and strongly in} \ L^2(\Omega),$$
	$$\mathrm{curl}\, \sigma(y_\alpha)\longrightarrow 
	\omega_0 \qquad
	 \mbox{weakly in} \ L^2(\Omega). $$
By taking into account $(\ref{var_form_state})$ and $(\ref{transport_state})$, we have
	\begin{equation}\label{vf_yalpha}\nu\left(\nabla y_\alpha,\nabla\phi\right)+\left(\mathbf{curl}\, \sigma(\mathbf{y_\alpha})\times \mathbf{y_\alpha},\boldsymbol\phi\right)=
	\left(u_\alpha,\phi\right) \qquad \mbox{for all} \ 
	\phi\in V\end{equation}
and 
	\begin{equation}\label{vf_sigma_alpha}\left(\mathrm{curl}\, \sigma(y_\alpha),\phi\right)-
	\tfrac{\alpha}{\nu}\, b(y_\alpha,\phi,\mathrm{curl}\, \sigma(y_\alpha))=\left(\tfrac{\alpha}{\nu}\,\mathrm{curl}\, u_\alpha+
	\mathrm{curl}\, y_\alpha,\phi\right)\qquad 
	\mbox{for all} \ 
	\phi\in {\cal D}(\Omega).\end{equation}
The previous convergence results yield
	$$\lim_{\alpha\rightarrow 0^+}\left(\mathbf{curl}\, \sigma(\mathbf{y_\alpha})\times \mathbf{y_\alpha},\boldsymbol\phi\right)=\left(\boldsymbol{\omega}_0\times \mathbf{y}_0,\boldsymbol\phi\right)\qquad \mbox{for all} \ 
	\phi\in V$$
and 
	$$\lim_{\alpha\rightarrow 0^+} 
	b(y_\alpha,\phi,\mathrm{curl}\, \sigma(y_\alpha))=
	b(y_0,\phi, \omega_0)\qquad 
	\mbox{for all} \ 
	\phi\in {\cal D}(\Omega).$$
Therefore, by passing to the limit in the previous identities, we obtain
	$$\nu\left(\nabla y_0,\nabla \phi\right)+\left(\boldsymbol{\omega}_0\times
	 \mathbf{y}_0,\boldsymbol\phi\right)=
	\left(u_0,\phi\right) \qquad \mbox{for all} \ 
	\phi\in V$$
and 
	$$\left(\omega_0,\phi\right)=\left(\mathrm{curl}\, 
	y_0,\phi\right)\qquad 
	\mbox{for all} \ 
	\phi\in {\cal D}(\Omega)$$
showing that $\omega_0=\mathrm{curl}\, 
	y_0$ and that $y_0$ satisfies
	$$\nu\left(\nabla y_0,\nabla \phi\right)
	+b(y_0,y_0,\phi)=
	\left(u_0,\phi\right) \qquad \mbox{for all} \ 
	\phi\in V.$$ 
that is, $(u_0,y_0)$ is admissible for $(P_0)$. Let us now prove that the convergence of $y_\alpha$ to $y_0$ is strong. Taking into account the variational formulations corresponding to $y_\alpha$ and $y_0$, we easily see that that $y_\alpha-y_0$ satisfies
	$$\begin{array}{ll}\nu \left|y_\alpha-y_0\right|_{H^1}^2
	&=\left(u_\alpha-u_0,y_\alpha- y_0\right)-
	\left(\mathbf{curl}\, \sigma(\mathbf y_\alpha)\times \mathbf y_0
	-\mathrm{curl}\, \mathbf y_0\times \mathbf y_0,
	\mathbf y_\alpha-\mathbf y_0\right)\vspace{2mm}\\
	&\longrightarrow 0 \qquad \mbox{when} \ \alpha \rightarrow 0^+.
	\end{array}$$
%%%%%%%%%%%%%%%%%%%%%%%%%%%%%%%%%%%%%%%%%%%%%%%%%
{\it Step 2. Convergence to an optimal solution of $(P_0)$.} Let us now prove that the limit point $(\bar u_0,\bar y_0)$ of a solution $(\bar u_\alpha,\bar y_\alpha)$ of $(P_\alpha)$ is a solution of $(P_0)$. By taking into account the convergence results established in the first step and the lower semicontinuity of $J$, we obtain
	$$\min(P_0)\leq J(\bar u_0,\bar y_0)\leq \liminf_{\alpha\rightarrow 0^+}J(\bar u_\alpha,\bar y_\alpha)=\liminf_{\alpha\rightarrow 0^+}\min(P_\alpha).$$ 
On the other hand, let $(\hat u,\hat y)$ be a solution of problem $(P_0)$ and let $\hat y_\alpha$ be the solution of $(\ref{equation_etat})$ corresponding to $\hat u$. Then $(\hat u,\hat y_\alpha)$ is admissible for $(P_\alpha)$ and 
	\begin{equation}\label{semi_con_1}
	\min(P_\alpha)\leq J(\hat u,\hat y_\alpha).\end{equation}
Arguing as in the first step, we can establish the convergence of $\hat y_\alpha$ to
 $\hat y$ in $V$ and thus 
	\begin{equation}\label{semi_con_2}\lim_{\alpha\rightarrow 0^+}\min(P_\alpha)\leq \lim_{\alpha\rightarrow 0^+}J(\hat u,\hat y_\alpha)=J(\hat u,\hat y)=\min(P_0).\end{equation}
Combining (\ref{semi_con_1}) and (\ref{semi_con_2}), we deduce that
	\begin{equation}\label{stability_alpha}
	\lim_{\alpha\rightarrow 0^+}\min(P_\alpha)=\min(P_0).\end{equation}
To prove the strong convergence of $\bar u_\alpha$ to $\bar u_0$ in $L^2(\Omega)$, observe that
	$$\begin{array}{ll}
	\left\|\bar u_\alpha-\bar u_0\right\|_2^2&=
	\left\|\bar u_\alpha\right\|_2^2
	-\left\|\bar u_0\right\|_2^2-
	2\left(\bar u_\alpha-\bar u_0,\bar u_0\right)
	\vspace{2mm}\\
	&=\tfrac{2}{\lambda}\left(
	J(\bar u_\alpha,\bar y_\alpha)-
	J(\bar u_0,\bar y_0)\right)-\tfrac{1}{\lambda}\left(
	\left\|\bar y_\alpha-y_d\right\|_2^2-
	\left\|\bar y_0-y_d\right\|_2^2\right)
	-2\left(\bar u_\alpha-\bar u_0,\bar u_0\right)\vspace{2mm}\\
	&=\tfrac{2}{\lambda}\left(
	\min(P_\alpha)-\min(P_0)\right)-\tfrac{1}{\lambda}\left(
	\left\|\bar y_\alpha-y_d\right\|_2^2-
	\left\|\bar y_0-y_d\right\|_2^2\right)
	-2\left(\bar u_\alpha-\bar u_0,\bar u_0\right).\end{array}$$
Therefore, by taking into account the convergence results of Step 1 and (\ref{stability_alpha}), it follows that
	\begin{equation}\label{L2_con_alpha}\lim_{\alpha\rightarrow 0^+}
	\left\|\bar u_\alpha-\bar u_0\right\|_2=0.\end{equation}
%%%%%%%%%%%%%%%%%%%%%%%%%%%%%%%%%%%%%%%%%%%%%%%%%
{\it Step 3. Convergence of $\bar p_\alpha$.} 
Due to (\ref{est_H1_p}), we have
	$$\left(1-\tfrac{ \bar\kappa}{\nu^2}\left(
	\left\|\bar u_\alpha\right\|_2+\alpha\left\|{\rm curl}\,\bar u_\alpha\right\|_2\right)
	\right) \left|\bar p_\alpha\right|_{H^1}\leq
	 \tfrac{\kappa}{\nu}\left\|\bar y_\alpha-y_d\right\|_2,$$
and by taking into account (\ref{L2_con_alpha}), we deduce that
 $\left(\bar p_\alpha\right)_\alpha$ is also bounded independently of $\alpha$. There then exists a subsequence, still indexed by $\alpha$ and $\bar p_0\in V$ such that
	$$\bar p_\alpha \longrightarrow \bar{p}_0 \qquad
	 \mbox{weakly in} \ H^1(\Omega) \
	\mbox{and strongly in} \ L^2(\Omega).$$
By taking into account the convergence results established in the first step, we deduce that
	$$\begin{array}{ll}\displaystyle\lim_{\alpha\rightarrow 0^+}\left(b\left( \bar p_\alpha,\phi, \sigma(\bar y_\alpha) \right)-b\left( \phi,\bar p_\alpha, \sigma(\bar y_\alpha) \right)\right)&=\displaystyle\lim_{\alpha\rightarrow 0^+}\left(\mathbf{curl}\,\sigma(\bar{\mathbf y}_\alpha) \times \bar{\mathbf p}_\alpha,\boldsymbol\phi\right)\vspace{2mm}\\
	&=\displaystyle \left(\mathbf{curl}\,\bar {\mathbf y}_0\times \bar{\mathbf p}_0,\boldsymbol\phi\right)\vspace{2mm}\\
	&=b(\bar p_0,\bar y_0,\phi)-b(\phi,\bar y_0,\bar p_0)
		\end{array}$$
and 
	$$\lim_{\alpha\rightarrow 0^+}\left(b\left(\bar p_\alpha, \bar y_\alpha, \sigma(\phi) \right)-b\left( \bar y_\alpha,\bar p_\alpha, \sigma(\phi) \right)\right)=b\left(\bar p_0, \bar y_0,\phi \right)-b\left(\bar  y_0,\bar p_0,\phi\right)$$
for all $\phi\in V$. Passing then to the limit in the variational formulation corresponding to $\bar p_\alpha$ yields
$$
	\nu\left(\nabla\bar p_0,\nabla\phi\right)+b(\phi,\bar y_0,\bar p_0)-b\left(\bar y_0,\bar p_0,\phi\right)
	=\left(\bar y_0-y_d, \phi\right)$$
for all $\phi\in V$ and thus $\bar p_0$ is the unique weak solution of 
	(\ref{adjoint_limit}). The optimality condition for the control follows then by passing to the limit in $(\ref{opt_control_alpha})$.$\hfill\Box$

\vspace{5mm}

\noindent {\bf Acknowledgment.} This work was partially supported by FCT project PEst-OE/MAT/UI4032/2011.

%%%%%%%%%%%%%%%%%%%%%%%%%%%%%%%%%%%%%%%%%%%%%%%%%%%%%%%%%%%%%%%%%%%


\smallskip

\begin{thebibliography}{99}

\bibitem{A12} {\sc N. Arada}, {\em Optimal control of shear-thinning fluids}, SIAM J. Control Optim. 50 (2012), pp. 2515-2542.

\bibitem{A13} {\sc  N. Arada}, {\em Optimal control of shear-thickening flows}, SIAM J. Control Optim. 51 (2013), 
pp. 1940-1961.

\bibitem{A14} {\sc  N. Arada}, {\em Optimal control of evolutionary quasi-Newtonian fluids}, SIAM J. Control Optim. 52 (2014), 
pp. 3401-3436.

%\bibitem{BI06} {\sc A. V. Busuioc, D. Iftimie},  {\em A non-Newtonian fluid with Navier boundary conditions}, J. Dynam. Diff. Eq. 18 (2006), pp. 357-379.

%\bibitem{BILN12}  {\sc A. V. Busuioc, D. Iftimie, M. C. Lopes Filho, H. J. Nussenzveig Lopes}, {\em Incompressible Euler as a limit of complex fluid models with Navier boundary conditions}, J. Differential Equations 252 (2012), pp. 624-640.

%\bibitem{B99} A.V. Busuioc, On second grade fluids with vanishing viscosity, C. R. Acad. Sci. Paris I, 328 (1999) 1241-1246. 
%Also in Port. Math. (N.S.) 59 (1) (2002) 47–65.

\bibitem{B99}{\sc J. M. Bernard}, {\em Stationary problem of second-grade fluids in three dimensions: existence, uniqueness and regularity}, Math. Meth. Appl. Sci. 22 (1999), pp. 655-687.

\bibitem{BR03} {\sc A. V. Busuioc, T. S. Ratiu},  {\em The second grade fluid and averaged Euler equations with Navier-slip boundary conditions},
Nonlinearity 16 (2003), pp. 1119-1149.

\bibitem{CG97}  {\sc D. Cioranescu, V. Girault}, {\em Weak and classical solutions of a family of second grade fluids}, Int. J. Nonlinear Mech. 32 (1997), pp. 317-335.

\bibitem{CO84}  {\sc D. Cioranescu, E. H. Ouazar}, {\em Existence and uniqueness for fluids of second grade}, Nonlinear Partial Differential Equations and Their Applications (Coll\`ege de France Seminar, Paris, 1982/1983), 4 (Boston, MA: Pitman) (1984), pp. 178-197.


%\bibitem{CMR98}  {\sc T. Clopeau, A. Mikelic, R. Robert}, {\em On the vanishing viscosity limit for the $2D$ incompressible Navier-Stokes equations with the friction type boundary conditions}, Nonlinearity 11 (1998), pp. 1625-1636.

%\bibitem{C96}  {\sc J. M. Coron}, {\em On the controllability of the $2D$ incompressible Navier-Stokes equations with the Navier-slip
%boundary conditions}, ESAIM Control Optim. Calc. Var. 1 (1996), pp. 35-75.

%\bibitem{DV02}  {\sc L. Desvillettes, C. Villani}, {\em On a variant of Korn's inequality arising in statistical mechanics}, ESAIM Control
%Optim. Calc. Var. 8 (2002), 603-619.

\bibitem{DF74}  {\sc J. E. Dunn, R. L. Fosdick}, {\em Thermodynamics, stability and boundedness of fluids of complexity $2$ and fluids of second grade}, Arch. Rational Mech. Anal. 56 (1974), pp. 191-252.

\bibitem{DR95}  {\sc J. E. Dunn, K. R. Rajagopal}, {\em Fluids of differential type: Critical review and thermodynamical analysis}, Int J. Eng. Sci. 33 (1995), pp. 689-729.

\bibitem{FHT2}  {\sc C. Foias, D. D. Holm, E. S. Titi}, {\em The Navier-Stokes-$\alpha$ model of fluid turbulence}, Physica D 153 (2001), pp. 505-519.

\bibitem{FHT1}  {\sc C. Foias, D. D. Holm, E. S. Titi}, {\em The three dimensional viscous Camassa-Holm equations and their relation to the Navier-Stokes equations and the turbulence theory}, J. Dynam. Diff. Eq. 14 (2002), pp. 1-35.

\bibitem{GS94}  {\sc G. P. Galdi, A. Sequeira}, {\em  Further existence results for classical solutions of the equations of a second-grade fluid}, Arch. Rational Mech. Anal. 128 (1994), pp. 297-312.

\bibitem{GR86} {\sc V. Girault, P. A. Raviart}, {\em Finite element
methods for Navier-Stokes equations. Theory and Algorithms.}
Springer-Verlag, Berlin (1986).

\bibitem{GS99}  {\sc V. Girault, L. R. Scott}, {\em  Analysis of a two-dimensional grade-two fluid model with a tangential boundary condition}, J. Math. Pures Appl. 78 (1999), pp. 981-1011.

\bibitem{HMR981}  {\sc D. D. Holm, J. E. Marsden, R. S. Ratiu}, {\em Euler-Poincar\'e models of ideal fluids with nonlinear dissipation}, Phys. Rev. Lett. 349 (1998), pp. 4173-4177.

 \bibitem{HMR98}  {\sc D. D. Holm, J. E. Marsden, R. S. Ratiu}, {\em The Euler-Poincar\'e equations and semi-direct products with applications to continuum theories}, Adv. Math. 137 (1998), pp. 1-81. 

%\bibitem{I02} D. Iftimie, Remarques sur la limite $\alpha\rightarrow 0$ pour les fluides de grade 2, C. R. Acad. Sci. Paris I 334 (1) (2002) 83-86. 
%Also in Studies in Mathematics and its Applications, 31 Nonlinear partial Differential Equations and their Applications. College de France Seminar, Volume XIV. North Holland, 2002.

%\bibitem{IP06} {\sc D. Iftimie, G. Planas}, {\em Inviscid limits for the Navier-Stokes equations with Navier friction boundary conditions}, Nonlinearity 19 (2006), pp. 899-918.

%\bibitem{K06}  {\sc J. P. Kelliher}, {\em Navier-Stokes equations with Navier boundary conditions for a bounded domain in the plane}, SIAM J. Math. Anal. 38 (2006), pp. 210-232.

%\bibitem{LT10}  {\sc J. S. Linshiz, E. S. Titi}, {\em On the convergence rate of the Euler-$\alpha$, an inviscid second-grade complex fluid, model to Euler equations}, J. Stat. Phys. 138 (2010), pp. 305-332.

%\bibitem{LNP05}  {\sc M.C. Lopes Filho, H.J. Nussenzveig Lopes, G. Planas}, {\em On the inviscid limit for $2D$ incompressible flow with Navier friction condition}, SIAM J. Math. Anal. 36 (2005), pp. 1130-1141.

\bibitem{LNTZ15}  {\sc M.C. Lopes Filho, H. J. Nussenzveig Lopes, E. S. Titi, A. Zang}, {\em Convergence of the $2D$ Euler-$\alpha$ to Euler equations in the Dirichlet case: Indifference to boundary layers}, Physica D 292-293 (2015) pp. 51-61. 

\bibitem{NT65}  {\sc W. Noll, C. Truesdell}, {\em The Nonlinear Field 
Theories of Mechanics}, Encyclopedia of Physics, (ed. S. Flugge), Vol. III/3, Springer-Verlag. 1965.

\bibitem{O81}  {\sc E. H. Ouazar}, {\em Sur les Fluides de Second Grade}, Th\`ese 3\`eme Cycle, Universit\'e Pierre et Marie Curie, 1981.

\bibitem{RE55}  {\sc R. S. Rivlin, J. L. Ericksen}, {\em Stress-deformation relations for isotropic materials}, Arch. Rational Mech. Anal. 4 (1955), pp. 323-425.
  
%\bibitem{S73} {\sc V. E. \v{S}\v{c}adilov, V. A. Solonikov}, {\em On a boundary value problem for a stationary system of Navier-Stokes equations}, Proc. Steklov Inst. Math. 125 (1973), pp. 186-199.

\bibitem{S05} {\sc T. Slawig}, {\em Distributed control for a class of non-Newtonian fluids}, J. Differential Equations 219 (2005), pp. 116-143.

\bibitem{WR10} {\sc D. Wachsmuth and T. Roub\'i\v{c}ek}, {\em Optimal control of incompressible non-Newtonian fluids}, Z. Anal. Anwend 29 (2010), pp. 351-376.


%\bibitem{S66}
%V. A. Solonikov, General boundary value problems for Douglis-Nirenberg elliptic systems, Proc. Steklov Inst. Math., \textbf{92} (1966), 269-339.






%\bibitem{JM01}
%Jager, W., Mikelic, A.: \emph{On the roughness-induced
%effective boundary conditions for a viscous flow.} J. Differential Equations,
%\textbf{170} (2001), 96--122.






%\bibitem[Lions et al.(1972)]{L72} Lions, P-L., Magenes, E.: \emph{Non-Homogeneous Boundary Value %Problems and Applications} Vol. \textbf{1}. Springer-Verlag . Heidelberg 
%(1972).





%\bibitem{Yu95}
%Yudovich, V. I.: \emph{Uniqueness theorem for the basic nonstationary problem in the dynamics of %an ideal incompressible fluid.}, Math. Res. Lett., \textbf{2} (1995), 27-38.

\end{thebibliography}
\end{document}